\documentclass[twocolumn]{autart}    

\usepackage{graphicx}          

\usepackage{amsmath,amssymb,bm,bbm,mathrsfs,amscd}
\usepackage{color}
\usepackage{dsfont}
\usepackage{graphicx}
\usepackage{epsfig}
\usepackage{subfigure}
\usepackage[noend]{algpseudocode}
\usepackage{tikz}
\usepackage{enumitem}
\usepackage{stackrel}
\usepackage{hyperref}
\usepackage{balance}
\newtheorem{theorem}{Theorem}
\newtheorem{lemma}{Lemma}
\newtheorem{remark}{Remark}
\newtheorem{example}{Example}
\newcommand{\RR}{{\mathbb{R}}}
\newcommand{\NN}{{\mathbb{N}}}
\newcommand{\EE}{{\mathbb{E}}}
\newcommand{\PP}{{\mathbb{P}}}
\newcommand{\JJ}{{\mathbb{J}}}
\newcommand{\FF}{{\mathbb{F}}}

\newcommand{\mc}{\mathcal}
\newcommand{\norm}[1]{\|#1\|}
\newcommand{\normsq}[1]{\|#1\|^2}
\newcommand{\EEk}[1]{\EE[#1|\mc F_k]}
\newcommand{\EEx}[1]{\EE[#1]}
\newcommand{\op}{\operatorname}

\newcommand{\prox}{\op{prox}}

\newcommand{\bs}{\boldsymbol}

\newcommand{\normsqphi}[1]{\|#1\|_\Phi^2}

\DeclareSymbolFont{myletters}{OML}{ztmcm}{m}{it}
\DeclareMathSymbol{\uplambda}{\mathord}{myletters}{"15}

\begin{document}

\begin{frontmatter}

\title{Stochastic generalized Nash equilibrium seeking under partial-decision information\thanksref{footnoteinfo}} 

\thanks[footnoteinfo]{This paper was not presented at any IFAC 
meeting. Corresponding author B.~Franci. Tel. +31-152785019. }

\author[Delft]{Barbara Franci}\ead{b.franci-1@tudelft.nl} and    
\author[Delft]{Sergio Grammatico}\ead{s.grammatico@tudelft.nl}               

\address[Delft]{Delft Center for Systems and Control, Delft University of Technology, Delft, The Netherlands}  


\begin{abstract}                          
We consider for the first time a stochastic generalized Nash equilibrium problem, i.e., with expected-value cost functions and joint feasibility constraints, under partial-decision information, meaning that the agents communicate only with some trusted neighbors. 
We propose several distributed algorithms for network games and aggregative games that we show being special instances of a preconditioned forward-backward splitting method. We prove that the algorithms converge to a generalized Nash equilibrium when the forward operator is restricted cocoercive by using the stochastic approximation scheme with variance reduction to estimate the expected value of the pseudogradient.
\end{abstract}

\end{frontmatter}

\section{Introduction}
In a stochastic Nash equilibrium problem (SNEP), some agents interact with the aim of minimizing their expected value cost function which is affected by the decision variables of the other agents. The characteristic feature is the presence of uncertainty, represented by a random variable with unknown distribution. Due to this complication, the equilibrium problem is typically hard to solve \cite{koshal2013,ravat2011}. However, many practical problems must be modelled with uncertainty, for instance, electricity markets with unknown demand \cite{henrion2007} and transportation systems with erratic travel time \cite{watling2006}.\\
Another recurrent engineering aspect is that agents may be subject to shared feasibility constraints. For instance, consider the gas market where the companies participate in a bounded capacity market \cite{abada2013} or more generally any network Cournot game with market capacity constraints and uncertainty in the demand \cite{yu2017,yi2019}. In this case, we have a stochastic generalized NEP (SGNEP), i.e., the problem of finding a Nash equilibrium where not only the cost function but also the feasible set depend on the decisions of the other agents \cite{chen2020,franci19,yi2019}. \\
This class of problems is of high interest for the decision and control community, in both deterministic \cite{chen2020,pavel2019,yi2019,gadjov2020} and stochastic cases \cite{koshal2013,yousefian2012}. Notably, the presence of shared constraints makes the computation of an equilibrium very challenging, especially when searching for distributed algorithms, where each agent only knows its local cost function and its local constraints. 
Perhaps the most elegant way to design a solution algorithm for a SGNEP is to recast the problem as an inclusion, leveraging on monotone operator theory. In particular, operator splitting methods, paired with a primal-dual analysis on the pseudogradient mapping, can be used to obtain fixed-point iterations that converge to an equilibrium, i.e., a collective strategy that simultaneously solve the interdependent optimization problems of the agents while reaching consensus on the dual variables associated to the coupling constraints \cite{facchinei2010,kulkarni2012}. Unfortunately, these methods do not necessarily lead to distributed iterations, thus, for this purpose, a suitable ``preconditioning'' is required in some problem classes \cite{belgioioso2018,yi2019}. \\
Although computationally convenient, distributed algorithms have one main flaw: the information that each agent must know about the others.
Nevertheless, in most cases, it is assumed that the agents have direct access to the decision variables of the other agents that affect their cost function \cite{franci19,koshal2013,ravat2011}. This is the so-called full-decision information setup, where the agents must share information with all the competing agents. From a more realistic point of view however, it is more natural to assume that the agents agree to share information with some trusted neighbors only. In this case, we have the so-called partial-decision information setup \cite{galeotti2010,pavel2019}. 
In the literature of deterministic GNEPs, there are several algorithms for both the full information setup \cite{facchinei2010,grammatico2017,yi2019} and the partial-decision information one \cite{belgioioso2018,bianchi2020,pavel2019}. On the contrary, to the best of our knowledge, the only few algorithms for SGNEPs are in full information \cite{franci19,ravat2011,yu2017}.
Among others, one of the fastest and less computationally demanding algorithms that may be exploited is the forward-backward (FB) splitting method \cite[Section 26.5]{bau2011}, for which a suitable preconditioning is needed to obtain distributed iterations \cite{franci19,yi2019}. In the stochastic case, the FB algorithm converges \cite{franci19,rosasco2016,yousefian2012} when the operator used for the forward step is strongly monotone \cite{franci19,rosasco2016} or cocoercive as we preliminarily show in the full decision information setup \cite{franci2020}.\\
There is another important issue in SGNEPs due to the presence of shared constraints: the monotonicity properties of the involved mappings are not necessarily preserved in the extended primal-dual operators obtained to decouple the shared constraints. Hence, ensuring convergence can be difficult because of the lack of a strongly monotone forward operator, not even when the original pseudogradient mapping is strongly monotone. Instead, cocoercivity can be obtained from a strongly monotone or cocoercive pseudogradient \cite{franci19,franci2020,yi2019}.
Nonetheless, in partial-decision information, the extended operator can only be at most restricted cocoercive with respect to the solution set and only when the pseudogradient mapping is strongly monotone \cite{gadjov2020,pavel2019}.\\
Besides the monotonicity properties, another challenging aspect is the uncertainty. The presence of a random variable with unknown distribution implies that the agents cannot compute the exact cost function but only an approximation \cite{franci2020,franci19,iusem2017,koshal2013}. This results in a stochastic error which complicates the analysis and prevents from applying the proofing techniques used in the deterministic case \cite{yi2019,pavel2019,gadjov2020,yi2020}.\\
In this paper we propose the first distributed algorithms specifically tailored for SGNEPs in partial-decision information and show their convergence to an equilibrium under restricted cocoercivity of the stochastic forward operator. Our contributions are summarized next:
\begin{itemize}
\item We model and study for the first time SGNEPs under partial-decision information. 
\item We propose two distributed algorithms for network games and two for aggregative games. The algorithms are characterized by the specific way we impose consensus on the dual variables, i.e., node-based or edge-based. While both the approaches are present in the literature of deterministic GNEPs \cite{bianchi2020,pavel2019,yi2020} they have been only partially used in the stochastic case and only in full information \cite{franci19,franci2020}.
\item We show that our algorithms are instances of a preconditioned FB splitting and we prove their convergence when the forward operator is \textit{restricted} cocoercive with respect to the solution set. The restricted cocoercivity assumption is much weaker than the monotonicity assumptions usually adopted in the stochastic literature \cite{rosasco2016,franci2020}.
\end{itemize}
As a special case, we also consider aggregative games, where the cost function of each agent does not depend explicitly on the individual decision of the other agents but it is related instead to some aggregate value of all the decisions \cite{gadjov2020,grammatico2017,koshal2016}. Illustrative examples are traffic networks where the time delay depends on the overall congestion \cite{paccagnan2019} and energy markets where the price of electricity depends on the aggregate demand \cite{chen2014}. In this case as well, the literature is not extensive \cite{lei2020,lei2018,meigs2019}.
While the authors in \cite{lei2020,lei2018} propose a stochastic proximal gradient response for aggregative games, relatively similar to ours, they do not consider shared constraints. Moreover, they assume a strongly monotone mapping to prove convergence. However, when dealing with SGNEPs in partial information, also in the aggregative case, the extended operator can be at most restricted cocoercive; therefore, the results in \cite{lei2020,lei2018} are not applicable to our generalized setting.\\
We remark that, despite our proposed algorithms are all instances of a FB scheme, hence inspired by the literature on the topic \cite{franci2020,franci19,pavel2019,gadjov2020}, this is the first time SGNEPs in partial decision information are addressed from an algorithmic point of view. Moreover, the edge-based approach is loosely studied even in the deterministic case, while here we show that it can be a valid alternative to the more classic node-based algorithm.\\
{\bf{Paper organization.}} The next section recalls some preliminary notions on operator and graph theory. SGNEPs in partial-decision information are described in Section \ref{sec_sgnep}. The first two algorithms for network games are presented in Section \ref{sec_p} while the aggregative case is discussed in Section \ref{sec_agg}. Sections \ref{sec_FB}, \ref{sec_conv_p} and \ref{sec_conv_a} are devoted to the theoretical convergence results. Specifically, in Section \ref{sec_FB} a fundamental lemma is proven and then, it is used in Sections \ref{sec_conv_p} and \ref{sec_conv_a} to show that the algorithms for network games and aggregative games, respectively, converge to an equilibrium. Numerical simulations (Section \ref{sec_sim}) and conclusion (Section \ref{sec_conclu}) end the paper.

\section{Preliminaries and Notation}
We use the same notation as in \cite{franci2020} and, with a slight abuse of notation, given two sets $A$ and $B$, we may indicate the Cartesian product as $\left[\begin{smallmatrix}A\\B\end{smallmatrix}\right]=A\times B$, to ease the reading.
The definitions are taken from \cite{bau2011,facchinei2007,godsil2013}.\\
{\bf{Monotone operator theory.}}
A mapping $F:\op{dom}F\subseteq\RR^n\to\RR^n$ is $\ell$-\textit{Lipschitz continuous} if, for some $\ell>0$, $\norm{F(x)-F(y)} \leq \ell\norm{x-y} \text { for all } x, y \in \op{dom}(F)$; $\eta$-\textit{strongly monotone} if, for some $\eta>0$, $\langle F(x)-F(y),x-y\rangle \geq \eta\|x-y\|^{2} \text{ for all }x, y \in \op{dom}(F);$
$\beta$-\textit{cocoercive} if, for some $\beta>0$, $ \langle F(x)-F(y),x-y\rangle \geq \beta\|F(x)-F(y)\|^{2}$, for all $x, y \in \op{dom}(F)$;
\textit{maximally monotone} if there exists no monotone operator $G :C\to\RR^n$ such that the graph of $G$ properly contains the graph of $F$ \cite[Def. 20.20]{bau2011}.
We use the adjective \textit{restricted} if a property holds for all $(x,y)\in\op{dom}(F)\times \op{fix}(F)$.\\
{\bf{Graph theory.}}
Basic definitions can be found in \cite{godsil2013}.
A graph $\mc G=(\mc I,\mc E)$ is undirected if $(i,j)\in\mc E$ and $(j,i)\in\mc E$. It is connected if there is a path between every pair of vertices. Let W be the weighted adjacency matrix, $L$ be the associated Laplacian matrix and $V$ be the incidence matrix. Then, if the graph is undirected $W$ and $L$ are symmetric, i.e., $W=W^\top$ and $L=L^\top$. If the graph is also connected, it holds that $L = V^\top V$ and that 
\begin{equation}\label{eq_null}
\op{null}(V ) = \op{null}(L) = \{\kappa \bs 1_N : \kappa  \in \RR\},
\end{equation}
namely, the null space of $V$ and $L$ is the consensus subspace. The laplacian matrix $L$ has an eigenvalue equal zero and let the other eigenvalues ordered as $0<\uplambda_2(L)\leq\dots\leq\uplambda_N(L)$ where $d_{\text{max}}\leq\uplambda_N(L)\leq2d_{\text{max}}$ and $d_{\text{max}}=\max_{i\in\mc I}\{d_i\}$. Given $d_{\text{max}}$, it follows from the Baillon-Haddad Theorem that the Laplacian is $\frac{1}{2d_{\text{max}}}$-cocoercive.

\section{Stochastic generalized Nash equilibrium problems under partial-decision information}\label{sec_sgnep}
\subsection{Problem setup}
We consider a stochastic generalized Nash equilibrium problem (SGNEP), i.e., the problem of finding a Nash equilibrium when the cost functions are expected value functions and the agents, indexed by $\mc I=\{1,\dots,N\}$, are subject to coupling constraints. \\
Each agent $i\in\mc I$ has a decision variable $x_i\in\Omega_i\subseteq\RR^{n_i}$ and a local cost function defined as 
\begin{equation}\label{eq_J}
\JJ_i(x_i,\bs{x}_{-i}):=\EE_\xi[f_i(x_i,\bs{x}_{-i},\xi_i(\omega))] + g_{i}(x_{i}),
\end{equation}
for some measurable function $f_i:\mc \RR^{n}\times \RR^d\to \RR$ and $n=\sum_{i\in\mc I} n_i$. 
Each agent $i$ aims at minimizing its local cost function within its feasible strategy set $\Omega_i$.
The cost function is split in smooth ($f_i$) and non smooth parts ($g_{i}:\RR^{n_i}\to\bar\RR$). The latter may also model local constraints via an indicator function ($g_{i}(x_{i})=\iota_{\Omega_i}(x_i)$).
\begin{assum}\label{ass_G}(Local cost functions)
For each $i\in\mc I$, the function $g_i$ in \eqref{eq_J} is lower semicontinuous and convex and $\op{dom}(g_{i})=\Omega_i$ is nonempty, compact and convex.
\end{assum}
For each agent $i$, the cost function in \eqref{eq_J} depends on the local variable $x_i$, on the decision of the other agents $\bs x_{-i}=\op{col}((x_j)_{j\neq i})$ and on the random variable $\xi_i(\omega)\in\RR^d$\footnote{From now on, we use $\xi$ instead of $\xi(\omega)$ and $\EE$ instead of $\EE_\xi$.}. 
The probability space is $(\Xi, \mc F, \PP)$ where $\Xi=\Xi_1\times \dots\times\Xi_N$.  We assume that the expected value $\EE[f_i(\bs{x},\xi)]$ is well defined for all the feasible $\bs{x}=\op{col}((x_i)_{i\in\mc I})\in\bs\Omega\subseteq\RR^n$, where $\bs\Omega=\prod_{i\in\mc I}\Omega_i$.\\
Besides the local constraints $x_i\in\Omega_i\subseteq\RR^{n_i}$, the agents are also subject to coupling constraints $A\bs x\leq b$, therefore, the feasible decision set of each agent $i \in \mc I$ is 
\begin{equation}\label{eq_X_i}
\mc{X}_i(\bs{x}_{-i}) :=\textstyle{\{y_i \in \Omega_i\; |\; A_i y_i \leq b-\sum_{j \neq i}^{N} A_j x_j\},}
\end{equation}
where $A_i \in \RR^{m \times n}$ and $b\in\RR^m$, and 
the collective feasible set reads as
\begin{equation}\label{eq_X}
\bs{\mc{X}}=\left\{\bs y \in\bs\Omega\; | \;A\bs y-b \leq \bs 0_{m}\right\},
\end{equation}
where $A=\left[A_{1}, \ldots, A_{N}\right]\in\RR^{m\times n}$, $A_i \in \RR^{m \times n_i}$ and $b\in\RR^m$.
We suppose that the constraints are deterministic and satisfy the following assumption \cite{facchinei2010,ravat2011}.
\begin{assum}\label{ass_X}(Constraint qualification)
For each $i \in \mc I,$ the set $\Omega_{i}$ is nonempty, compact and convex.
The set $\bs{\mc{X}}$ satisfies Slater's constraint qualification. 
\end{assum}
Given the decision variables of the other agents $\bs{x}_{-i}$, the goal of each agent $i$ is to choose a strategy $x_i$ that solves its local optimization problem, i.e.,
\begin{equation}\label{eq_game}
\forall i \in \mc I: \quad \left\{\begin{array}{cl}
\min_{x_i \in \Omega_i} & \JJ_i\left(x_i, \bs{x}_{-i}\right) \\ 
\text { s.t. } & A_i x_i \leq b-\sum_{j \neq i}^{N} A_j x_j.
\end{array}\right.
\end{equation}
By simultaneously solving all the coupled optimization problems, we have a stochastic generalized Nash equilibrium (SGNE).
\begin{defn}
A stochastic generalized Nash equilibrium is a collective strategy $\bs x^*\in\bs{\mc X}$ such that for all $i \in \mc I$
$$\JJ_i(x_i^{*}, \bs x_{-i}^{*}) \leq \inf \{\JJ_i(y, \bs x_{-i}^{*})\; | \; y \in \mc{X}_i(\bs x_{-i}^{*})\}.$$
\end{defn}
To guarantee the existence of a SGNE \cite[Section 3.1]{ravat2011}, we make further assumptions on the cost function, typical of the deterministic setup as well \cite{ravat2011,facchinei2007vi,facchinei2010}.
\begin{assum}\label{ass_J}(Cost functions convexity)
For every $i \in \mc I$ and $\bs{x}_{-i} \in \bs{\mc{X}}_{-i}$ the function $f_{i}(\cdot, \bs{x}_{-i})$ is convex and continuously differentiable.
For every $i\in\mc I$ and for every $\xi_i \in \Xi_i$, the function $f_{i}(\cdot,\bs x_{-i},\xi_i)$ is convex, continuously differentiable, and Lipschitz continuous and for each $\bs x_{-i}$; the Lipschitz constant $\ell_i(\bs x_{-i},\xi_i)$ is integrable in $\xi_i$. The function $f_{i}(x_i,\bs x_{-i},\cdot)$ is measurable.
\end{assum}

Among all the possible equilibria, we focus on the class named variational equilibria (v-SGNE), i.e., those equilibria that are also solution of a suitable stochastic variational inequality (SVI). To describe this class, let us introduce the pseudogradient mapping
\begin{equation}\label{eq_F}
\FF(\bs{x})=\op{col}\left((\EE[\nabla_{x_{i}} f_{i}(x_{i}, \bs{x}_{-i},\xi)])_{i\in\mc I}\right),
\end{equation} 
where the exchange between the expected value and the gradient is possible because of Assumption \ref{ass_J} \cite{ravat2011}.
A standard assumption on the pseudogradient in partial-decision information is postulated next \cite{gadjov2020,pavel2019}.
\begin{assum}\label{ass_strong}(Strongly monotone pseudogradient)
$\FF$ in \eqref{eq_F} is $\eta$-strongly monotone and $\ell_\FF$-Lipschitz continuous, for some constants $\eta,\ell_\FF>0$ respectively.
\end{assum}
\begin{example}
An affine map $F(x)=Ax+b$ with $A\in\RR^{n\times n}$ symmetric and positive definite is strongly monotone \cite[Page 155]{facchinei2007}.
\end{example}

\begin{remark}\label{remark_sol}
Under Assumption \ref{ass_strong}, the associated SVI has a unique solution, therefore, there is a unique v-SGNE \cite[Theorem 2.3.3]{facchinei2007}, \cite[Lemma 3.3]{ravat2011}.
\end{remark}

As in \cite{facchinei2010,auslender2000,ravat2011}, the SGNEP can be recasted as the monotone inclusion 
\begin{equation}\label{eq_T}
0\in\mc T(\bs{x},\bs\lambda):=\left[\begin{smallmatrix}
G(\bs{x})+\FF(\bs{x})+A^{\top} \lambda \\ 
\mathrm{N}_{\RR_{ \geq 0}^{m}}(\lambda)-(A \bs{x}-b)
\end{smallmatrix}\right],
\end{equation}
i.e., as the problem of finding a zero of the set-valued mapping $\mc T:\mc X\times \RR^m_{\geq 0}\rightrightarrows \RR^{n}\times\RR^m$, where $G(\bs{x})=\partial g_{1}(x_{1})\times \ldots\times \partial g_{N}(x_{N})$. The operator in \eqref{eq_T} can be obtained via a primal-dual characterization of the equilibria: the $i$-th component of the first row of $\mc T$ in \eqref{eq_T} corresponds to $0 \in \EE[\nabla_{x_i} f_i(x_i^{*}, \bs{x}_{-i}^{*},\xi_i)]+\partial g_i(x_i^{*})+A_i^{\top} \lambda$, for $i\in\mc I$, i.e., to the stationarity condition of each optimization problem in \eqref{eq_game} while the second row is the complementarity condition. Indeed, a collective decision $\bs x^{*}$ is a v-SGNE of the game in \eqref{eq_game} if and only if the Karush--Kuhn--Tucker (KKT) conditions associated to \eqref{eq_game} are satisfied with consensus of the dual variables, i.e., $\lambda_i=\lambda$ for all $i\in\mc I$ \cite[Theorem 3.1]{facchinei2007vi}, \cite[Theorem 3.1]{auslender2000}.\\
Moreover, we consider a partial-decision information setup where the agents have access only to some of the other players decision variables, exchanged locally over an undirected communication graph $\mc G= (\mc I , \mc E )$. 
\begin{assum}\label{ass_graph}(Graph connectivity)
The graph $\mc G=(\mc I,\mc E)$ is undirected and connected.
\end{assum}
To overcome the lack of knowledge of the decision that affects its cost, each agent keeps an estimate of the action of the other players \cite{gadjov2020,pavel2019}. Let us denote with $\hat x_{i,j}$ the estimate of the decision of agent $j$ stored by agent $i$ and let us collect all the estimates stored by agent $i$ in $\hat {\bs x}_i=\op{col}((\hat x_{i,j})_{j\in\mc I})\in\RR^{n}$. We note that $\hat x_{i,i}=x_i$ and let $\hat {\bs x}_{i,-i}=\op{col}((\hat x_{i,j})_{j\neq i})$.
Thus, to compute the equilibria of the game in \eqref{eq_game}, the agents should reach consensus not only on the dual variables, but also on the estimates, i.e., $\lambda_i=\lambda_j$ and $\hat{\bs x}_i=\hat{\bs x}_j$ for all $i,j\in\mc I$. With this aim, let us introduce the consensus subspace of dimension $q\in\NN$
$$\vspace{-.3cm}C_q=\{\bs y\in \RR^{Nq} :\bs y=1_N\otimes y,y\in \RR^q\}$$
and its orthogonal complement $C_q^\perp$. Then, $C_n$ is the consensus subspace of the estimated decisions while $C_m$ is the consensus subspace of the dual variables.\\
To compute the v-SGNE of the game in \eqref{eq_game}, the agents use the estimates, i.e., the pseudogradient mapping is modified according to
\begin{equation}\label{eq_F_p}
\begin{aligned}
\FF_p(\hat{\bs{x}})&=\op{col}((\nabla_{x_{i}} \JJ_{i}(x_{i}, \hat{\bs{x}}_{i,-i}))_{i \in \mathcal{I}})\\
&=\op{col}((\EE[\nabla_{x_{i}}f_i(x_{i}, \hat{\bs{x}}_{i,-i},\xi_i)])_{i \in \mathcal{I}}).
\end{aligned}
\end{equation}\vspace{-.5cm}
\begin{remark}\label{remark_lip_p}
The extended pseudogradient $\FF_p$ in \eqref{eq_F_p} is $\ell_p$-Lipschitz continuous with constant $0<\ell_p\leq\ell_\FF$, as a consequence of Assumption \ref{ass_strong} \cite[Lemma 3]{bianchi2020ecc}.
\end{remark}
\subsection{Approximation scheme}
Since the random variables have an unknown distribution, i.e., the expected values in \eqref{eq_F} are virtually impossible to compute, we take an approximation of the pseudogradient mapping.
We suppose that the agents have access to an increasing number $M_k$ of samples of the random variables $\xi_i$ and that they are able to compute an approximation of $\FF(\bs x)$ of the form $F^\textup{SA}(\bs x,\bs\xi)=\op{col}(F_i^{\textup{SA}}(\bs x,\bar\xi_i)),$ where 
\begin{equation}\label{eq_F_SA}
\textstyle{F_i^\textup{SA}(\bs x,\bar\xi_i)=\frac{1}{M_k} \sum_{t=1}^{M_k} \nabla_{x_i}f_i(\bs x,\xi_i^{(t)}),}
\end{equation}
$\bar\xi_i=\op{col}(\xi_i^{(1)},\dots,\xi_i^{(M_k)})$ for all $i\in\mc I$ and $\bs \xi=\op{col}(\bar \xi_1,\dots,\bar\xi_n)$ is an i.i.d. sequence of random variables drawn from $\PP$.
Approximations of the form (\ref{eq_F_SA}), using a finite (increasing) number of samples, are known as stochastic approximations (SA) and they are very common in Monte-Carlo simulation, machine learning and computational statistics \cite{iusem2017}.
From now on, we indicate with the superscript SA the operators where the mapping $\FF$ is sampled with $F^\textup{SA}$ as in \eqref{eq_F_SA}. 
\begin{assum}\label{ass_batch}(Increasing batch size)
The batch size sequence $(M_k)_{k\geq 1}$ is such that, for some $c,k_0,a>0$,
$M_k\geq c(k+k_0)^{a+1}.$
\end{assum}
From the last assumption it follows that $1/M_k$ is summable, which is standard when using a SA scheme, especially in combination with the forthcoming variance reduction assumption (Assumption \ref{ass_error}) \cite{iusem2017}. \\
Since we use an approximation, for $k \geq 0,$ let us introduce the stochastic error 
\begin{equation}\label{eq_stoc_error}
\epsilon_k=F^\textup{SA}(\bs x^k,\bs \xi^k)-\FF(\bs x^k).
\end{equation}
Let us define the filtration $\mc F=\{\mc F_k\}_{k\in\NN}$, that is, a family of $\sigma$-algebras such that $\mathcal{F}_{0} = \sigma\left(X_{0}\right)$, for all $k \geq 1$, $\mathcal{F}_{k} = \sigma\left(X_{0}, \xi_{1}, \xi_{2}, \ldots, \xi_{k}\right)$ and $\mc F_k\subseteq\mc F_{k+1}$ for all $k\geq0$.
Then, standard assumptions for the stochastic error are to have zero mean and bounded variance \cite{iusem2017,koshal2013}. 
\begin{assum}\label{ass_error}(Zero mean and bounded variance)
For all $k\geq 0$, for all $\bs x\in\bs{\mc X}$ 
$\EEk{\epsilon_k}=0,$ a.s.
and there exists $\sigma>0$ such that for all $\bs x\in\bs{\mc X}$
\begin{equation}\label{eq_variance}
\EEx{\norm{F^\textup{SA}(\bs x, \cdot)-\FF(\bs x)}^{2}}\leq \sigma^2.
\end{equation}
\end{assum}

\begin{remark}\label{rem_variance}
Under Assumptions \ref{ass_batch} and \ref{ass_error}, it holds that for all $k>0$
\begin{equation}\label{eq_vr}
\mathbb{E}[\|\epsilon_{k}\|^{2} \mid \mathcal{F}_{k}] \leq \tfrac{c \sigma^{2}}{M_{k}},
\end{equation}
where $M_k$ is the batch size sequence used in the approximation \eqref{eq_F_SA} (see \cite[Lem. 3.12]{iusem2017}, \cite[Lem. 6]{franci2020} for a proof). Since \eqref{eq_vr} implies that the second moment of the error diminishes with the number of samples $M_k$, algorithms using the approximation in \eqref{eq_F_SA} are also known as variance-reduced methods \cite{iusem2017}.
\end{remark}

\section{Stochastic preconditioned forward-backward algorithms for network games}\label{sec_p}
In this section, we present two distributed algorithms for network games.
We suppose that each agent $i\in\mc I$ only knows its own cost function $\JJ_i$, its feasible set $\Omega_i$, and its own portion of the coupling constraints $A_i$ and $b_i$.
Moreover, through the graph $\mc G$, the agents have access to some of the variables of the other agents.
In Section \ref{sec_conv_p} we show that the algorithms are instances of a preconditioned forward-backward (pFB) algorithm \cite{yi2019,bau2011} and we show how to choose suitable operators to derive them.
\subsection{Node-based algorithm for network games}\label{sec_pn}
We start with the distributed iterations presented in Algorithm \ref{algo_pn}. 
Its steps involve: a proximal step to update each decision variable $x_i$; an updating rule for the estimates that pushes $\hat{\bs x}_i$ toward consensus; the auxiliary variable $z_i$ which helps reaching the dual variable consensus \cite{yi2019}; a projection step into the positive orthant for the dual variable $\lambda_i$. 
\begin{alg}\label{algo_pn}
{{(Node-based fully-distributed preconditioned forward-backward)}}\\
\noindent\rule{8.5cm}{.4pt}\\
Initialization: $x_i^0 \in \Omega_i, \lambda_i^0 \in \RR_{\geq0}^{m},$ and $z_i^0 \in \RR^{m} .$\\
Iteration $k$: Agent $i$\\
(1) Receives $x_j^k$ and $\lambda_j^k$ for $j \in \mathcal{N}_{i}$, then updates
$$
\begin{aligned}
x_i^{k+1}=&\op{prox}_{g_i}[x_i^k-\alpha_{i}(F^\textup{SA}(x_i^k, \hat{\bs x}_{i,-i}^k,\bar\xi_i^k)+A_{i}^\top \lambda_i^k\\
&+c\textstyle{\sum_{j \in \mathcal{N}_{i}} w_{i j}(x_i^k-\hat{\bs x}_{i,j}^k))] }\\
\hat{\bs x}_{i, -i}^{k+1}=&\hat{\bs x}_{i,-i}^k-\alpha_{i}c \textstyle{\sum_{j \in \mathcal{N}_{i}} w_{i j}(\hat{\bs x}_{i,-i}^k-\hat{\bs x}_{j,-i}^k )}\\
z_i^{k+1}=&z_i^k-\nu_{i}\textstyle{ \sum_{j \in \mathcal{N}_{i}} w_{i j}(\lambda_i^k-\lambda_{j, k}) }
\end{aligned}
$$
(2) Receives $x_j^{k+1}$ and $z_j^{k+1}$ for $j \in \mathcal{N}_{i}$, then updates
$$
\begin{aligned}
\lambda_i^{k+1}&=\op{proj}_{\RR_{\geq0}^{m}}[\lambda_i^k+\delta_{i}(A_{i}(2 x_i^{k+1}-x_i^k)-b_{i} \\
&-\textstyle{\sum_{j \in \mathcal{N}_{i}} w_{i j}(2(z_i^{k+1}-z_{j}^{ k+1})-(z_i^k-z_{j}^{k})))]} \\
\end{aligned}\vspace{-.2cm}
$$
\noindent\rule{8.5cm}{.4pt}\\
\end{alg}\vspace{-.5cm}
The algorithm is inspired by the preconditioned FB iterations proposed in \cite{pavel2019}. The main difference is that Algorithm \ref{algo_pn} is not deterministic, thus, for the update of the primal variable, the approximation of the extended pseudogradient mapping in \eqref{eq_F_p} is used.
The algorithm is fully distributed since each agent $i$ knows its own variables and shares its information only with the neighbors in $\mc N_i$.
The algorithm is characterized by the choice of the consensus constraint for the dual variables. In this case, exploiting \eqref{eq_null}, it is imposed as $L\lambda=0$ where $L$ is the laplacian matrix associated to $\mc G$. Since we use the laplacian matrix, we call the algorithm \textit{node-based}. Similarly, it is also imposed the consensus constraint on the estimates, $\bs L\hat{\bs x}=0$, in the update of $\hat{\bs x}$.\\
We now state our first convergence result.
\begin{theorem}\label{theo_pn}
Let Assumptions \ref{ass_G}--\ref{ass_error} hold. Then, there exist $\bar \alpha, \bar \nu,\bar \delta>0$ such that, for $\alpha_i\in(0,\bar\alpha)$, $\nu_i\in(0,\bar\nu)$ and $\delta_i\in(0,\bar\delta)$, for all $i\in\mc I$, the sequence $(\bs x^k)_{k\in\NN}$ generated by Algorithm \ref{algo_pn} converges a.s. to a v-SGNE of the game in \eqref{eq_game}.
\end{theorem}\vspace{-.7cm}
\begin{pf} 
See Section \ref{sec_conv_pn} where we also provide explicit bounds for the step sizes.
\qed\end{pf}


\subsection{Edge-based algorithm for network games}\label{sec_pe}

Let us now describe another instance of the pFB algorithm that differs from Algorithm \ref{algo_pn} in the way we impose the consensus constraint on the dual variables. Specifically, following \eqref{eq_null}, we impose the constraint $V\lambda=0$. The details on how we exploit this \textit{edge-based} constraint, i.e., using the incidence matrix, are presented in Section \ref{sec_conv_pe} while the iterations are presented in Algorithm \ref{algo_pe}.
\begin{alg}\label{algo_pe}
{{(Edge-based fully-distributed preconditioned forward-backward)}}\\
\noindent\rule{8.5cm}{.4pt}\\
Initialization: $x_i^0 \in \Omega_i, \lambda_i^0 \in \RR_{\geq0}^{m},$ and $z_i^0 \in \RR^{m} .$\\
Iteration $k$: Agent $i$\\
(1) Receives $x_j^k$ and $\lambda_j^k$ for $j \in \mathcal{N}_{i}$, then updates
$$
\begin{aligned}
x_{i}^{k+1}=& \op{prox}_{g_{i}}[x_{i}^{k}-\alpha_{i}(F^\textup{SA}_{i}(x_{i}^{k},\hat{\bs x}_{i,-i}^{k},\bar\xi_i^k)+A_i^{\top} \lambda_{i}^{k}\\ 
&\textstyle{+c\sum_{j\in\mc N_i}w_{ij}(x_i^k-\hat{\bs x}_{i,j}^k))] }\\
\hat{\bs x}_{i, -i}^{k+1}=&\hat{\bs x}_{i,-i}^k-\alpha_{i} \textstyle{\sum_{j \in \mathcal{N}_{i}} w_{i j}(\hat{\bs x}_{i,-i}^k-\hat{\bs x}_{j,-i}^k )}\\
z_{i}^{k+1}=& z_{i}^{k}-\nu(\lambda_{i}^{k}-\textstyle{\sum_{j\in\mc N_i} w_{i j} \lambda_{j}^{k})} \\
\lambda_{i}^{k+1}=&\op{proj}_{\mathbb{R}_{\geq 0}^{m}}[\lambda_{i}^{k}+\delta_{i}(2 A_i x_{i}^{k+1}-A_i x_{i}^{k}-b_i\\
&-2 z_{i}^{k+1}+z_{i}^{k})]
\end{aligned}\vspace{-.2cm}
$$
\noindent\rule{8.5cm}{.4pt}\\
\end{alg}\vspace{-.5cm}\vspace{-.2cm}
A consequence of the edge-based constraint is that only one communication round is required because each $\lambda_i$ depends only on local variables. The updating rule of $x_i$ and $\hat{\bs x}_i$, instead, are the same as in Algorithm \ref{algo_pn} because they are not affected by the dual variable constraint.\\
The use of the incidence matrix is not common, even in the deterministic case. Similar iterations has been considered in \cite{yi2020} which however proposes a deterministic asynchronous Krasnoselskii-Mann iteration (see also Section \ref{sec_conv_pe}).\\
We can state the convergence result for Algorithm \ref{algo_pe}.
\begin{theorem}\label{theo_pe}
Let Assumptions \ref{ass_G}--\ref{ass_error} hold. Then, there exist $\bar \alpha, \bar \nu,\bar \delta>0$ such that, for $\alpha_i\in(0,\bar\alpha)$, $\nu_i\in(0,\bar\nu)$ and $\delta_i\in(0,\bar\delta)$, for all $i\in\mc I$, the sequence $(\bs x^k)_{k\in\NN}$ generated by Algorithm \ref{algo_pe} converges a.s. to a v-SGNE of the game in \eqref{eq_game}.
\end{theorem}\vspace{-.7cm}
\begin{pf}
See Section \ref{sec_conv_pe} where we also provide explicit bounds for the step sizes.
\qed\end{pf}

\section{Stochastic aggregative games}\label{sec_agg}

With aggregative games we mean a class of games where the cost function explicitly depends on the aggregate decision of all the agents. 
Formally, given the actions of all the players $x_i\in\RR^{n_i}$ where $n_i=\bar n$ for all $i\in\mc I$, let 
$$\textstyle{\op{avg}(\bs x)=\frac{1}{N}\sum_{i\in\mc I}x_i}$$
be the average strategy.  Then, the cost function of each agent $i\in\mc I$ in the aggregative case can be written as
$$\JJ_i(x_i,x_{-i})=\EE[f_i(x_i,\op{avg}(\bs x),\xi_i)]+g_i(x_i)$$
where $f_i:\RR^n\times\RR^d\to\RR$ satisfies Assumption \ref{ass_J}, $g_i$ is as in Assumption \ref{ass_G} and $\xi_i:\Xi_i\to\RR^d$ is the uncertainty. Notice that in this case, since $n_i=\bar n$ for all $i\in\mc I$, $\op{avg}(\bs x)\in\RR^{\bar n}$ and $n=\bar nN$.\\
Since this is a particular case of the classical SGNEP in \eqref{eq_game}, existence and uniqueness of an equilibrium hold under the same assumptions and the v-SGNE can be characterized using the KKT conditions in \eqref{eq_T}.
Accordingly, Algorithm \ref{algo_pn} and Algorithm \ref{algo_pe} can be used to reach an equilibrium. However, the previous algorithms require the agents to exchange the estimates of all the other actions, i.e., a vector of dimension $(N-1)\bar n$, while the aggregate value has dimension $\bar n$ (independent of the number of agents). To reduce the computational complexity, we propose two algorithms, depending on the consensus constraint, tailored for aggregative games.\\
Let us introduce the pseudogradient mapping for the aggregative case as $\FF_a(\bs x,y)=\op{col}(\FF^a_i(x_i,y_i)_{i\in\mc I}),$
where
\begin{equation}\label{eq_F_a}
\textstyle{\FF^a_i(x_i,y_i)=\EE[\nabla_{x_{i}} f_{i}(x_{i},y_i,\xi)+\frac{1}{N}\nabla_{y_{i}} f_{i}(x_{i},y_i,\xi)].}
\end{equation}
The variable $y$ indicates the dependency on the aggregate value. In fact, $\FF_i(x_i,\op{avg}(\bs x))=\nabla_{x_{i}}\EE[f_i(x_i,\op{avg}(\bs x),\xi)]=\nabla_{x_{i}}\JJ_i(x_i,x_{-i})$, i.e., $\FF_a(\bs x,\op{avg}(\bs x))=\FF(\bs x)$.
\begin{remark}\label{remark_lip_a}
It follows from Assumptions \ref{ass_strong} and Remark \ref{remark_lip_p} that
$\FF_a$ in \eqref{eq_F_a} is Lipschitz continuous in both the arguments with constants $\ell^x_a,\ell^u_a>0$, respectively.
\end{remark}
Due to the partial-decision information setup, the agents cannot compute the exact average strategy.
To overcome this problem, each agent updates an auxiliary variable $s_i=\op{avg}(\bs x)-x_i\in\RR^{\bar n}$ \cite{bianchi2020}.
The variable $\bs s=\op{col}(s_i)_{i\in\mc I}$ is used to reconstruct the true aggregate value, controlling only the information received from the neighboring agents. Specifically, it should hold that $\bs s_k\to\bs 1_N\otimes \op{avg}(\bs x_k )-\bs x_k$ asymptotically. Moreover, let
\begin{equation}\label{eq_u}
u_{i}:=x_{i}+s_{i}, 
\end{equation}
and $\bs u:=\op{col}((u_{i})_{i \in \mathcal{I}})$.
The variable $\bs u^k$ represents the approximated average through the iterations. We remark that the explicit tracking of $\bs u^k$ is not necessary in our algorithms since, to estimate di aggregative value, we update iteratively the variable $\bs s^k$. 
Moreover, from the updating rule of $\bs s^k$, it follows that $\op{avg}(\bs s^k)=\bs 0_{\bar n}$. Thus, provided that the algorithm is initialized appropriately, i.e., $s_{i}^{0}=\mathbf{0}_{\bar{n}},$ for all $i \in \mathcal{I}$, an invariance property holds for the approximated average, i.e., for all $k\in\NN$,
\begin{equation}\label{eq_invariance}
\op{avg}(\bs x^k)=\op{avg}(\bs u^k).\vspace{-.5cm}
\end{equation}


\subsection{Node-based algorithm for aggregative games}\label{sec_an}
We first consider the node-based consensus constraint introduced in Section \ref{sec_pn}. 
Since in this case we have to take into consideration also the aggregative value, the state variable is $\bs{\omega}=\op{col}(\bs x, \bs s,\bs z, \bs \lambda)$, where $\bs x$ is the exact decision variable, $\bs s$ is the tracking variable, $\bs z$ is  the auxiliary variable for consensus of the dual variables and $\bs\lambda$ is the dual variable.

\begin{alg}\label{algo_an}
{{(Node-based fully-distributed preconditioned forward-backward for aggregative games)}}\\
\noindent\rule{8.5cm}{.4pt}\\
Initialization: $x_i^0 \in \Omega_i, \lambda_i^0 \in \RR_{\geq0}^{m},$ and $z_i^0 \in \RR^{m} .$\\
Iteration $k$: Agent $i$\\
(1) Receives $x_j^k$, $s_j^k$ and $\lambda_j^k$ for $j \in \mathcal{N}_{i}$, then updates
$$\begin{aligned}
x_i^{k+1}=&\op{prox}_{g_i}[x_i^k-\alpha_{i}(F^\textup{SA}_{i}(x_i^k, x_{i}^k+s_i^k,\bar\xi_i^k)+A_{i}^\top \lambda_i^k\\
&\textstyle{+\sum_{j \in \mathcal{N}_{i}} w_{i,j}(x_{i}^k+s_i^k-(x_{j}^k+s_j^k)))]}\\
s_i^{k+1}=&s_i^k-\gamma_i \textstyle{\sum_{j \in \mathcal{N}_{i}} w_{i,j}(x_{i}^k+s_i^k-(x_{j}^k+s_j^k))}\\
z_i^{k+1}=&z_i^k-\nu_{i} \textstyle{\sum_{j \in \mathcal{N}_{i}} w_{i,j}(\lambda_i^k-\lambda_{j}^k)}\\
\end{aligned}$$
(2) Receives $z_j^{k+1}$ for $j \in \mathcal{N}_{i}$, then updates
$$\begin{aligned}
\lambda_i^{k+1}=&\op{proj}_{\RR_{+}^{m}}\left[\lambda_i^k+\delta_{i}\left(A_{i}(2x_i^{k+1}-x_i^k)-b_{i}\right)\right.+\\
&-\delta_{i}\textstyle{\sum_{j \in \mathcal{N}_{i}} w_{i,j}\left(2(z_i^{k+1}-z_{j}^{k+1})-(z_i^k-z_{j}^k)\right)]}\\
\end{aligned}\vspace{-.2cm}$$
\noindent\rule{8.5cm}{.4pt}
\end{alg}\vspace{-.4cm}
Compared to Algorithm \ref{algo_pn}, besides the presence of the variable $\bs s$, Algorithm \ref{algo_an} has a different updating rule for $\bs x^k$, that now includes the estimated aggregative value $\bs x^k+\bs s^k$. The remaining variables ($\bs z^k$ and $\bs \lambda^k$) are not influenced by the average strategy, therefore the updating rules are the same as in Algorithm \ref{algo_pn}. 
The operators used to obtain the iterations in Algorithm \ref{algo_an} are given in Section \ref{sec_conv_an} but here we state the convergence result.
\begin{theorem}\label{theo_an}
Let Assumptions \ref{ass_G}--\ref{ass_error} hold. Then, there exist $\bar \alpha, \bar \nu,\bar \delta>0$ and $\gamma_i>0$ such that, for $\alpha_i\in(0,\bar\alpha)$, $\nu_i\in(0,\bar\nu)$ and $\delta_i\in(0,\bar\delta)$, for all $i\in\mc I$, the sequence $(\bs x^k)_{k\in\NN}$ generated by Algorithm \ref{algo_an} converges a.s. to a v-SGNE of the game in \eqref{eq_game}.
\end{theorem}\vspace{-.7cm}
\begin{pf}
See Section \ref{sec_conv_an} where we also provide explicit bounds for the step sizes. \qed
\end{pf}
\begin{remark}\label{remark_gadjov}
In \cite{gadjov2020}, the authors propose an algorithm in which the agents keep track of the whole aggregative value $\bs u$ (instead of $\bs s$). 
Specifically, the aggregative value is updated as
\begin{equation}\label{eq_u_k}
\bs u^{k+1}=\bs u^{k}-\gamma \bs{L}_{\bar{n}} \bs u^{k}+(\bs x^{k+1}-\bs x^{k})
\end{equation}
where $\bs{L}_{\bar{n}}:=L \otimes I_{\bar{n}}$. The rule in \eqref{eq_u_k} can be derived from the updating rule of $\bs s^k$ in Algorithm \ref{algo_an} or \ref{algo_ae} and its definition in \eqref{eq_u} and it can be regarded as a dynamic tracking of the time-varying quantity $\op{avg}(\bs x)$ \cite{gadjov2020,koshal2016}.
The algorithm in \cite{gadjov2020} is still an instance of a pFB but the operators and preconditioning matrix are different from ours. 
See Section \ref{sec_conv_a} for further technical details.
\end{remark}


\subsection{Edge-based algorithm for aggregative games}\label{sec_ae}
In this section, we consider an edge-based pFB algorithm, similarly to Section \ref{sec_pe}. The iterations in the aggregative case read as in Algorithm \ref{algo_ae}.
\begin{alg}\label{algo_ae}
{{(Edge-based fully-distributed preconditioned forward-backward for aggregative games)}}\\
\noindent\rule{8.5cm}{.4pt}\\
Initialization: $x_i^0 \in \Omega_i, \lambda_i^0 \in \RR_{\geq0}^{m},$ and $z_i^0=0$.\\
Iteration $k$: Agent $i$\\
(1) Receives $x_j^k$, $s_j^k$ and $\lambda_j^k$ for $j \in \mathcal{N}_{i}$, then updates
$$
\begin{aligned}
x_{i}^{k+1}=& \op{prox}_{g_{i}}[x_{i}^{k}-\alpha_{i}(F^\textup{SA}_{i}(x_{i}^{k}, x_{i}^k+s_i^k,\bar \xi_i^k)+A_i^{\top} \lambda_{i}^{k} \\
&\textstyle{+c\sum_{j \in \mathcal{N}_{i}} w_{i,j}(x_{i}^k+s_i^k-(x_{j}^k+s_j^k)))]}\\
s_i^{k+1}=&s_i^k-\gamma_i \textstyle{\sum_{j \in \mathcal{N}_{i}} w_{i,j}(x_{i}^k+s_i^k-(x_{j}^k+s_j^k))}\\
z_{i}^{k+1}=& z_{i}^{k}-\nu(\lambda_{i}^{k}-\textstyle{\sum_{j\in\mc N_i} w_{i j} \lambda_{j}^{k}) }\\
\lambda_{i}^{k+1}=&\op{proj}_{\mathbb{R}_{\geq 0}^{m}}[\lambda_{i}^{k}+\delta_{i}(2 A_i x_{i}^{k+1}-A_i x_{i}^{k}-b_i\\
&-2 z_{i}^{k+1}+z_{i}^{k})]
\end{aligned}\vspace{-.2cm}
$$
\vspace{-.7cm}
\noindent\rule{8.5cm}{.4pt}\\
\end{alg}
We note that the updating rule for $x_i$ and for the auxiliary variable $s_i$ are the same as in Algorithm \ref{algo_an} while the difference is in the auxiliary variable $z_i$ and the dual variable $\lambda_i$ that now depends only on local variables. This follows using the edge-based consensus constraint as in Section \ref{sec_pe}. 
More details on how to obtain the iterations and the proof of the following result can be found in Section \ref{sec_conv_ae}.

\begin{theorem}\label{theo_ae}
Let Assumptions \ref{ass_G}--\ref{ass_error} hold. Then, there exist $\bar \alpha, \bar \nu,\bar \delta>0$ such that, for $\alpha_i\in(0,\bar\alpha)$, $\nu_i\in(0,\bar\nu)$ and $\delta_i\in(0,\bar\delta)$, $\gamma_i>0$ for all $i\in\mc I$, the sequence $(\bs x^k)_{k\in\NN}$ generated by Algorithm \ref{algo_ae} converges a.s. to a v-SGNE of the game in \eqref{eq_game}.
\end{theorem}\vspace{-.7cm}
\begin{pf}
See Section \ref{sec_conv_ae} where we also provide explicit bounds for the step sizes.
\qed\end{pf}

\section{Convergence analysis: A fundamental lemma}\label{sec_FB}

In this section we show that the classic pFB splitting converges a.s. in the stochastic case to the zeros of the operator $\mc T$ in \eqref{eq_T} when the forward operator is restricted cocoercive. First of all, let us rewrite the operator $\mc T$ into the summation of the two operators 
\begin{equation}\label{eq_op}
\begin{aligned}
\mc{A} &:\left[\begin{smallmatrix}
\bs{x} \\
\bs\lambda
\end{smallmatrix}\right] \mapsto\left[\begin{smallmatrix}
\FF(\bs{x}) \\ 
b
\end{smallmatrix}\right],\\
\mc{B} &:\left[\begin{smallmatrix}
\bs{x} \\ 
\bs\lambda
\end{smallmatrix}\right] \mapsto\left[\begin{smallmatrix}
G(\bs{x}) \\ 
\mathrm{N}_{\RR_{ \geq 0}^{m}}(\lambda)
\end{smallmatrix}\right]+\left[\begin{smallmatrix}
0 & A^{\top} \\ 
-A & 0
\end{smallmatrix}\right]\left[\begin{smallmatrix}
\bs{x} \\ 
\lambda
\end{smallmatrix}\right].
\end{aligned}
\end{equation}
Then, finding a solution of the SGNEP translates in finding a pair $(\bs x^*,\bs \lambda^*)\in\bs{\mc X}\times \RR_{\geq0}^m$ such that $(\bs x^*,\bs\lambda^*)\in\op{zer}(\mc A+\mc B)$. 
The zeros of the mapping $\mc T=\mc A+\mc B$ can be obtained through a FB splitting \cite[Section 26.5]{bau2011},\cite{belgioioso2018,yi2019},
that for any matrix $\Phi\succ0$, leads to the FB iteration:
\begin{equation}\label{eq_FB}
\bs\omega^{k+1}=(\mathrm{Id}+\Phi^{-1} \mc B)^{-1} \circ(\mathrm{Id}-\Phi^{-1} \mc A)(\bs\omega^k),
\end{equation}
where $\bs \omega^k=\op{col}(\bs x^k,\bs \lambda^k)$, $(\mathrm{Id}+\Phi^{-1} \mc B)^{-1}$ is the backward step and $(\mathrm{Id}-\Phi^{-1} \mc A)$ is the forward step.\\
We note that the convergence of \eqref{eq_FB} is independent on the specific choice of the operators $\mc A$ and $\mc B$ as long as some monotonicity conditions are satisfied. For this reason, we postulate the following assumption.

\begin{assum}\label{ass_coco}
The forward operator is restricted $\beta$-cocoercive for some $\beta>0$ and the backward operator is maximally monotone.
\end{assum}
\begin{remark}
An affine map $F(x)=Ax+b$ with $A\in\RR^{n\times n}$ symmetric and positive semidefinite is cocoercive \cite[Page 79]{facchinei2007}.
More generally, every $\eta$-strongly monotone, $\ell$-Lipschitz continuous function is $\frac{\eta}{\ell^2}$-cocoercive, but the vice-versa is not true in general.
\end{remark}

Checking if Assumption \ref{ass_coco} holds for suitable operators will be the key to prove convergence of the algorithms presented in the previous sections. Before stating the convergence result of the pFB iteration, however we need to postulate some further assumptions and to consider the approximation scheme. In fact, since the random variable has an unknown distribution, we replace $\mc A$ with $\mc A^\textup{SA}$, the operator obtained using the approximation in \eqref{eq_F_SA}. Thus, the pFB iteration reads as 
\begin{equation}\label{eq_FB_SA}
\bs\omega^{k+1}=(\mathrm{Id}+\Phi^{-1} \mc B)^{-1} \circ(\mathrm{Id}-\Phi^{-1} \mc A^\textup{SA})(\bs\omega^k).
\end{equation} 
We note that there is no uncertainty in the constraints, therefore, the corresponding error of the approximated extended operator is
$$\varepsilon_k=\mc A^\textup{SA}(\bs\omega^k,\bs\xi^k)-\mc A(\bs\omega^k)=\op{col}(\epsilon_k,0).$$
To guarantee convergence, the preconditioning matrix $\Phi$ should be positive definite. Since this property may depend on the specific choice of the matrix $\Phi$, here we postulate it as an assumption and in the following sections we ensure that it holds for the proposed algorithms.

\begin{assum}\label{ass_phi}
$\Phi$ is positive definite, i.e., $\Phi\succ0$.
\end{assum}
Moreover, to guarantee convergence and independently on the choice of $\Phi$, the step sizes should be bounded.
\begin{assum}\label{standass_phi}
$\norm{\Phi^{-1}}<2\beta$ where $\beta$ is the cocoercivity constant of the forward operator as in Assumption \ref{ass_coco}.
\end{assum}


We can now state and prove the convergence result for the iteration in \eqref{eq_FB_SA}.
\begin{lemma}\label{lemma_FB}
Let Assumptions \ref{ass_G} - \ref{standass_phi} hold. Then, the sequence $(\bs x_k,\bs \lambda_k)_{k\in\NN}$ generated by \eqref{eq_FB_SA} converges a.s. to some $(\bs x^*,\bs \lambda^*)\in\op{zer}(\mc A,\mc B)$ where $\bs x^*$ is a v-SGNE the game in \eqref{eq_game}, $\mc A$ and $\mc B$ are as in \eqref{eq_op}, and $\mc A^\textup{SA}$ is approximated using \eqref{eq_F_SA}.
\end{lemma}\vspace{-.7cm}
\begin{pf}
For brevity, we let $\hat{\mc A}=\mc A^\textup{SA}$.
We start by using that the resolvent is firmly nonexpansive \cite[Cor. 23.9]{bau2011} and that if $\bs\omega^*$ is a solution then it is a fixed point of the FB iteration in \eqref{eq_FB}: 
\begin{equation*}
\begin{aligned}
&\normsqphi{\bs\omega^{k+1}-\bs\omega^*}\leq
\normsqphi{\bs\omega^k-\bs\omega^*}+2\langle \bs\omega^k-\bs\omega^*,\Phi^{-1}\varepsilon_k\rangle_\Phi\\
&-2\langle \bs\omega^k-\bs\omega^*,\Phi^{-1}(\mc A(\bs\omega^k)-\mc A(\bs\omega^*))\rangle_\Phi+\\
&-\normsqphi{\bs\omega^k-\bs\omega^{k+1}}+2\langle\bs\omega^k-\bs\omega^{k+1},\Phi^{-1}(\hat{\mc A}(\bs\omega^k)-\mc A(\bs\omega^*))\rangle_\Phi.
\end{aligned}
\end{equation*}
Choosing $\zeta>1$ and such that Assumption \ref{ass_phi} is satisfied, we can use Young's inequality to obtain
\begin{equation}\label{young}
\begin{aligned}
&2\langle\bs\omega^k-\bs\omega^{k+1},\Phi^{-1}(\hat{\mc A}(\bs\omega^k)-\mc A(\bs\omega^{k+1}))\rangle_\Phi\leq\tfrac{1}{\zeta}\normsqphi{\bs\omega^k-\bs\omega^*}+\\
&+\zeta\normsqphi{\Phi^{-1}(\mc A(\bs\omega^k)-\mc A(\bs\omega^*))}+\zeta\normsqphi{\Phi^{-1}\varepsilon_k}+\\
&+2\zeta\langle \Phi^{-1}\mc A(\bs\omega^k)-\Phi^{-1}\mc A(\bs\omega^*),\varepsilon_k\rangle_\Phi.
\end{aligned}
\end{equation}
Then, by using restricted cocoercivity of $\mc A$ and including \eqref{young}, we obtain: 
\begin{equation}\label{step}
\begin{aligned}
&\normsqphi{\bs\omega^{k+1}-\bs\omega^*}
\leq\normsqphi{\bs\omega^k-\bs\omega^*}+\zeta\normsqphi{\Phi^{-1}\varepsilon^k}+\\
&+\textstyle{\left(\frac{1}{\zeta}-1\right)}\normsqphi{\bs\omega^k-\bs\omega^{k+1}}+2\langle \bs\omega^k-\bs\omega^*,\Phi^{-1}\varepsilon_k\rangle_\Phi+\\
&+\textstyle{\left(\frac{\zeta\norm{\Phi^{-1}}}{\theta}-2\right)}\langle \bs\omega^k-\bs\omega^*,\Phi^{-1}\mc A(\bs\omega^k)-\Phi^{-1}\mc A(\bs\omega^*)\rangle_\Phi\\
&+2\zeta\langle \Phi^{-1}\mc A(\bs\omega^k)-\Phi^{-1}\mc A(\bs\omega^*),\varepsilon_k\rangle_\Phi
\end{aligned}
\end{equation}
Next, let $\op{res}_\Phi(\bs\omega^k)^2=\normsqphi{\bs\omega^k-(\op{Id}+\Phi^{-1}\mc B)^{-1}(\bs\omega^k-\Phi^{-1}\mc A(\bs\omega^k))}$ and let us recall that $\op{res}_\Phi(\bs\omega)=0$ if and only if $\bs \omega$ is a solution \cite[Proposition 1.5.8]{facchinei2007}. It holds that
\begin{equation*}
\begin{aligned}
\op{res}_\Phi&(\bs\omega^k)^2=\normsqphi{\bs\omega^k-(\op{Id}+\Phi^{-1}\mc B)^{-1}(\bs\omega^k-\Phi^{-1}\mc A(\bs\omega^k))}\\
\leq&2\normsqphi{\bs\omega^k-\bs\omega^{k+1}}+\\
&+2\normsqphi{(\op{Id}+\Phi^{-1}\mc B)^{-1}(\bs\omega^k-\Phi^{-1}\hat{\mc A}(\bs\omega^k,\xi^k))+\\
&-(\op{Id}+\Phi^{-1}\mc B)^{-1}(\bs\omega^k-\Phi^{-1}\mc A(\bs\omega^k))}\\
\leq&
2\normsqphi{\bs\omega^k-\bs\omega^{k+1}}+2\normsqphi{\Phi^{-1}\varepsilon^k}
\end{aligned}
\end{equation*}
where the first equality follows by adding and subtracting $\bs \omega^{k+1}$ and using its definition and the last inequality follows from nonexpansivity. Then,
$\normsqphi{\bs\omega^k-\bs\omega^{k+1}}\geq\frac{1}{2}\op{res}_\Phi(\bs\omega^k)^2-\normsqphi{\Phi^{-1}\varepsilon^k}$. 
Finally, equation \eqref{step} becomes
\begin{equation*}
\begin{aligned}
&\normsqphi{\bs\omega^{k+1}-\bs\omega^*}\leq\normsqphi{\bs\omega^k-\bs\omega^*}+\textstyle{\left(\zeta-\frac{1}{\zeta}+1\right)}\normsqphi{\Phi^{-1}\varepsilon^k}+\\
&+2\langle \bs\omega^k-\bs\omega^*,\Phi^{-1}\varepsilon_k\rangle_\Phi+\textstyle{\frac{1}{2}\left(\frac{1}{\zeta}-1\right)}\op{res}_\Phi(\bs\omega^k)^2\\
&+\textstyle{\left(\frac{\zeta\norm{\Phi^{-1}}}{\theta}-2\right)}\langle \bs\omega^k-\bs\omega^*,\Phi^{-1}\mc A(\bs\omega^k)-\Phi^{-1}\mc A(\bs\omega^*)\rangle_\Phi\\
&+2\zeta\langle \Phi^{-1}\mc A(\bs\omega^k)-\Phi^{-1}\mc A(\bs\omega^*),\varepsilon_k\rangle_\Phi
\end{aligned}
\end{equation*}
By Assumption \ref{ass_phi} and by monotonicity, the second last term is smaller or equal than zero, hence, by taking the expected value and by using Assumption \ref{ass_error} and Remark \ref{rem_variance} we have that 
$$\begin{aligned}
\EEk{\normsqphi{\bs\omega^{k+1}-\bs\omega^*}}\leq&\normsqphi{\bs\omega^k-\bs\omega^*}+2\tfrac{c\sigma^2\norm{\Phi^{-1}}}{M_k}\\
&+\textstyle{\frac{1}{2}\left(\frac{1}{\zeta}-1\right)\op{res}(\bs\omega^k)^2}.
\end{aligned}$$
Using Robbins-Siegmund Lemma \cite{RS1971}
we conclude that $(\bs\omega^k)_{k\in\NN}$ is bounded and has a cluster point $\bar{\bs \omega}$. Moreover, it follows that $\sum_k\frac{1}{2}\left(\frac{1}{\zeta}-1\right)\op{res}_\Phi(\bs\omega^k)^2<\infty$, hence, $\op{res}(\bs\omega^k)\to0$ as $k\to\infty$ and $\op{res}(\bar{\bs \omega})=0$.
\qed\end{pf}
\begin{remark}
We note that the operators $\mc A$ and $\mc B$ is \eqref{eq_op} satisfy Assumption \ref{ass_coco}  \cite[Lemma 1]{belgioioso2018} and that the matrix
$$
\Phi=\left[\begin{smallmatrix}
\alpha^{-1} & -A^{\top} \\
-A & \gamma^{-1}
\end{smallmatrix}\right]
$$
is positive definite \cite[Lemma 3]{belgioioso2018}, therefore, the pFB algorithm in \eqref{eq_FB_SA} obtained with these operators converges to a v-SGNE of the game in \eqref{eq_game}. However, expanding \eqref{eq_FB_SA}, the iterations that we obtain require full information on the decision of the other agents \cite{belgioioso2018}. 
\end{remark}

\section{Convergence analysis for network games}\label{sec_conv_p}
We now show how to obtain suitable forward and backward operators that lead to Algorithms \ref{algo_pn} and \ref{algo_pe}, presented in Section \ref{sec_p}. Later, we show that such operators satisfy the assumptions of Lemma \ref{lemma_FB}, i.e., that the algorithms converge a.s. to a v-SGNE of the game in \eqref{eq_game}.\\
Let us first introduce some notation. Similarly to \cite{pavel2019}, let us define, for all $i\in\mc I$, the matrices
\begin{equation}\label{eq_rs}
\begin{aligned}
\mathcal{R}_{i} &:=\left[\begin{smallmatrix}
\mathbf{0}_{n_{i} \times n_{<i}} & I_{n_{i}} & \mathbf{0}_{n_{i} \times n_{>i}}
\end{smallmatrix}\right] \\
\mathcal{S}_{i}: &=\left[\begin{smallmatrix}
I_{n_{<i}} & \mathbf{0}_{n_{<i} \times n_{i}} & \mathbf{0}_{n_{<i} \times n_{>i}} \\
\mathbf{0}_{n_{>i} \times n_{<i}} & \mathbf{0}_{n_{>i} \times n_{i}} & I_{n_{>i}}
\end{smallmatrix}\right]
\end{aligned}
\end{equation}
where $n_{<i}:=\sum_{j<i, j \in \mathcal{I}} n_{j}, n_{>i}:=\sum_{j>i, j \in \mathcal{I}} n_{j} .$ The two matrices in \eqref{eq_rs} can be interpret as follows: $\mathcal{R}_{i}$ selects the $i$-th $n_{i}$ dimensional component from an $n$-dimensional vector, while $\mathcal{S}_{i}$ removes it. 
Namely, $\mathcal{R}_{i} \hat{\bs{x}}_{i}=\hat{\bs{x}}_{i, i}=x_{i}$ and $\mathcal{S}_{i} \hat{\bs{x}}_{i}=\hat{\bs{x}}_{i,-i} .$ Letting $\mathcal{R}:=\op{diag}((\mathcal{R}_{i})_{i \in \mathcal{I}})$ and $\mathcal{S}:=\op{diag}((\mathcal{S}_{i})_{i \in \mathcal{I}})$, we have that $\bs x=\mathcal{R} \hat{\bs{x}}$, $\op{col}((\hat{\bs{x}}_{i,-i})_{i \in \mathcal{I}})=\mathcal{S} \hat{\bs{x}}$ and $\hat{\bs{x}}=\mathcal{R}^{\top} \bs x+\mathcal{S}^{\top} \mathcal{S} \hat{\bs{x}}$ \cite{pavel2019}. \\
We can now analyze the two algorithms separately.

\subsection{Convergence of Algorithm \ref{algo_pn}}\label{sec_conv_pn}
Let $L\in\RR^{N\times N}$ be the Laplacian of the communication graph $\mc G$  and let $\bs{L}_m=L\otimes \op{Id}_m\in\RR^{Nm\times Nm}$ and $\bs{L}_n=L\otimes \op{Id}_n\in\RR^{Nn\times Nn}$. Let also ${\bf{A}}=\op{diag}\{A_1,\dots,A_N\}\in\RR^{Nm\times n}$ and $\bs \lambda=\op{col}(\lambda_1,\dots,\lambda_N)\in\RR^{Nm}$; similarly let $\bs b\in\RR^{Nm}$. 
As already mentioned in Section \ref{sec_pn}, exploiting \eqref{eq_null}, to impose consensus on the primal and dual variables, one can consider the Laplacian constraints $\bs{L}\hat{\bs x}=0$ and $\bs{L}\bs\lambda=0$ \cite{yi2019,pavel2019}. To include these constraints, similarly to \cite{pavel2019}, let us define the operators
\begin{equation}\label{eq_op_pn}
\begin{aligned}
\mc A_p: &\left[\begin{smallmatrix}
\hat{\bs x} \\
\bs z\\
\bs \lambda
\end{smallmatrix}\right]\hspace{-.1cm}\mapsto\hspace{-.1cm}\left[\begin{smallmatrix}
\mathcal{R}^\top \FF_p(\hat{\bs x})+c \bs L_{n} \hat{\bs x} \\
\mathbf{0} \\
\bs{b}
\end{smallmatrix}\right]\\
\mc B_p: &\left[\begin{smallmatrix}
\hat{\bs x} \\
\bs z\\
\bs \lambda
\end{smallmatrix}\right]\hspace{-.1cm}\mapsto\hspace{-.1cm}\left[\begin{smallmatrix}
\hat G(\bs{\hat x})  \\
\mathbf{0} \\
\op{N}_{\bs{R}_{+}^{N m}}(\bs{\lambda})
\end{smallmatrix}\right]\hspace{-.1cm}+\hspace{-.1cm}
\left[\begin{smallmatrix}
\mathcal{R}^\top {\bf{A}}^{\top}\lambda\\
-\bs L_m\lambda\\
-{\bf{A}} \mathcal{R}\hat{\bs x}+ \bs L_m \bs z
\end{smallmatrix}\right]
\end{aligned}
\end{equation}
where $\hat G(\hat{\bs x})=\mc R^\top G(\mc R\hat{\bs x})=\{\mc R^\top v:v\in G(\mc R\hat{\bs x})\}$ and $G(\mc R\hat{\bs x})=G(\bs x)=\partial g_1(x_1)\times\dots\times\partial g_N(x_N)$.
We note that we consider the estimates $\hat{\bs x}$ as a state variable and we use the matrix $\mc R$ to select the variables corresponding to each agent. Compared to the operators in \cite{pavel2019}, we have the expected valued extended pseudogradient in \eqref{eq_F_p} that, in the iterative process, is replaced by a stochastic approximation of the form in \eqref{eq_F_SA}.\\ 
The term $\bs L_n\hat{\bs x}$ is a measure of the disagreement between the decision variables of the agents and the estimates. Each agent can use this term to move towards consensus of the estimates while it uses the gradient to minimize the cost. We note that in a full information setup, the disagreement term can be removed (setting $c=0$) \cite{pavel2019}.\\
Given the operators in \eqref{eq_op_pn}, the pFB algorithm in compact form reads as
\begin{equation}\label{eq_FB_pn}
\bs\omega^{k+1}=(\mathrm{Id}+\Phi_p^{-1} \mc B_p)^{-1} \circ(\mathrm{Id}-\Phi_p^{-1} \mc A_p^\textup{SA})(\bs\omega^k),
\end{equation}
where the operator $\mc A_p^{\textup{SA}}$ is the operator $\mc A_p$ with the approximated pseudogradient mapping $F^\textup{SA}$ as in \eqref{eq_F_SA} and $\bs \omega^k=\op{col}(\hat{\bs x}^k,\bs z^k,\bs \lambda^k)$.
To obtain the distributed iterations in Algorithm \ref{algo_pn}, a suitable preconditioning matrix should be taken \cite{pavel2019}. Specifically, let
\begin{equation}\label{eq_phi_pn}
\Phi_p=\left[\begin{smallmatrix}
\alpha^{-1} & 0 & -\mathcal{R}^{\top} {\bf{A}}^{\top} \\
0 & \nu^{-1} & \bs L_m \\
-{\bf{A}} \mathcal{R} & \bs L_m & \delta^{-1}
\end{smallmatrix}\right]
\end{equation}
where $\alpha=\op{diag}\{\alpha_1\op{Id}_{n_1},\dots,\alpha_N\op{Id}_{n_N}\}\in\RR^{n\times n}$, and similarly $\nu$ and $\delta$ are block diagonal matrices collecting the step sizes. Then, expanding the iterations in \eqref{eq_FB_pn} with $\mc A_p$ and $\mc B_p$ as in \eqref{eq_op_pn} and $\Phi_p$ as in \eqref{eq_phi_pn}, we obtain 
\begin{equation}\label{compact_pn}
\begin{aligned}
\hat{\bs x}^{k+1}&\ni\hat{\bs x}^k\!-\!\alpha(\mc R^\top F^\textup{SA}(\hat{\bs x})+\mc R^\top {\bf{A}}^\top\bs\lambda^k+c\bs L_n\hat{\bs x})\!+\!\hat G(\hat{\bs x})\\
\bs z^{k+1}&=\bs z^k\!-\!\nu \bs L_m\bs \lambda^k\\
\bs\lambda^{k+1}&=\bs\lambda^k\!+\!\delta(A\mc R(2\hat{\bs x}^{k+1}\!-\hat{\bs x}^k)\!-\!b\!+\!\bs L_m(2\bs z^{k+1}-\bs z^k)).
\end{aligned}
\end{equation}
From the first line of \eqref{compact_pn} we obtain the update for both the decision variable $x_i^{k+1}$ of agent $i$ and the estimates $\hat{\bs x}_{i,-i}^{k+1}$ \cite[Lemma 1]{pavel2019}.
Precisely, premultiplying the first line of \eqref{compact_pn} by $\mc R$
we obtain $\bs x^{k+1}=\prox_G[\bs x^k-\alpha(F^\textup{SA}(\hat{\bs x})+ {\bf{A}}^\top\bs\lambda^k+c\mc R\bs L_{\bar n}\hat{\bs x})],$
that is, the update for each agent decision variable $x_i$ as in Algorithm \ref{algo_pn}. Instead, if we premultiply by $\mc S$
we obtain $\mc S\hat{\bs x}^{k+1}=\mc S\hat{\bs x}^k-\alpha c\mc S\bs L_{\bar n}\hat{\bs x},$
i.e., the update of the estimates.\\
For the sake of the convergence analysis, we have to guarantee that the preconditioning matrix $\Phi_p$ is positive definite (Assumption \ref{ass_phi}) and bounded in norm (Assumption \ref{standass_phi}), therefore, we take some bounds on the step sizes \cite{pavel2019}.

\begin{assum}\label{ass_phi_pn}
For a given $\tau>0$, the step sizes $\bar\alpha$, $\bar\nu$ and $\bar\delta$ are such that, for all $i\in\mc I$,
$$
\begin{aligned}
0&<\alpha_{i} \leq\bar\alpha\leq(\tau+\max\nolimits_{j\in\{1, \ldots n_{i}\}}\textstyle\sum_{k=1}^{m}|\left[A_{i}^{\top}\right]_{j k}|)^{-1} \\
0&<\nu_{i} \leq\bar\nu\leq(\tau+2 d_{i})^{-1}\\
0&<\delta_{i} \leq\bar\delta\leq(\tau+2 d_{i}+\max\nolimits_{j\in\{1, \ldots m\}}\textstyle\sum_{k=1}^{n_{i}}|\left[A_{i}\right]_{j k}|)^{-1} \\
\end{aligned}
$$
where $[A_i^\top]_{jk}$ indicates the entry $(j,k)$ of the matrix $A_i^\top$, and such that $\Phi_p$ satisfies Assumption \ref{standass_phi}.
\end{assum}
Then, it follows from the Gershgorin Theorem and \cite[Lemma 5]{pavel2019} that $\Phi_p\succ0$. 
We can now prove the convergence result. \vspace{-.6cm}
\begin{pf}[Proof of Theorem \ref{theo_pn}]
First, we show that the zeros of $\mc A_p+\mc B_p$ correspond to a v-SGNE of the game in \eqref{eq_game}. 
Expanding the inclusion $\bs \omega=\op{col}(\hat{\bs x}^*,\bs z^*,\bs\lambda^*)\in\op{zer}(\mc A_p+\mc B_p)$ we obtain
\begin{equation}\label{eq_incl_p}
\begin{aligned}
0 &\in \mathcal{R}^\top \FF_p(\hat{\bs x}^*)+c \bs L_n \hat{\bs x}^*+\mathcal{R}^\top {\bf{A}}^\top \boldsymbol{\lambda}^{*}+\hat G(\hat{\bs x}^*) \\
0&=-\bs L_m \bs\lambda^{*} \\
0 &\in \boldsymbol{b}+\op{N}_{\RR_{\geq0}^{N m}}(\boldsymbol{\lambda}^{*})-{\bf{A}} \mathcal{R} \hat{\bs x}^*+\bs L_m \boldsymbol{z}^{*}.
\end{aligned}
\end{equation}
Then, from the second line of \eqref{eq_incl_p} it follows that $\bs\lambda^*\in \op{null}(\bs L_m)$, i.e., $\bs\lambda^*=\mathbf{1}_N\otimes\lambda^*$ by \eqref{eq_null}. Similarly to \cite[Theorem 1]{pavel2019}, from the first line of \eqref{eq_incl_p} it follows that $\hat{\bs x}^*\in \op{null}(\bs L_n)$, i.e., $\hat{\bs x}^{*}=\mathbf{1}_{N} \otimes \bs x^{*}$ and that the first KKT condition in \eqref{eq_T} is satisfied. 
From the third line, we obtain the second KKT condition in \eqref{eq_T} \cite[Theorem 1]{pavel2019}.\\
Moreover, it also holds that $\op{zer}(\mc A_p+\mc B_p) \neq \varnothing$. In fact, by Remark \ref{remark_sol}, there exists a unique solution $\bs x^*$ and, therefore, there exists $\bs\lambda^*$ such that the KKT conditions in \eqref{eq_T} are satisfied \cite[Proposition 1.2.1]{facchinei2007} and $(\bs x^*,\bs\lambda^*)\in\op{zer}(\mc T)$. The first two lines of \eqref{eq_incl_p} are satisfied and using \eqref{eq_T} we can prove that there exists $\bs z$ such that the third line is satisfied as well \cite[Theorem 1]{pavel2019}.\\
We now show that the operators in \eqref{eq_op_pn} have the properties in Assumption \ref{ass_coco}. To this aim, we define some preliminary quantities \cite[Lemma 3]{pavel2019}. Let
$$
\Upsilon_p=\left[\begin{smallmatrix}
\frac{\eta}{N} & -\frac{\ell_p+\ell_\FF}{2 \sqrt{N}} \\
-\frac{\ell_p+\ell_\FF}{2 \sqrt{N}} &c\uplambda_{2}(L)-\ell_p
\end{smallmatrix}\right]
$$
and let $\mu_p=\uplambda_{\min}(\Upsilon_p)$ be its smaller eigenvalue. Let $c>c_{\min },$ where $c_{\min } \uplambda_{2}(L)=\frac{\left(\ell_p+\ell_\FF\right)^{2}}{4 \eta}+\ell_p$. Let us indicate with $Z=C_n\times\RR^{Nm}\times\RR^{Nm}$ the set where there is consensus on the first component, i.e., on the primal variable.
First, let us recall that from Assumption \ref{ass_phi_pn} it follows that $\Phi\succ0$ \cite[Lemma 5]{pavel2019}.
Then, the fact that $\mc A_p$ is $\beta_p$-restricted cocoercive with respect to $Z$, with constant $\beta_p\in(0,\frac{\mu_p}{\theta_p^2}]$, where $\mu_p=\uplambda_{\min}(\Upsilon_p)$ and $\theta_p=\ell_{p}+2cd_{\text{max}}$ follows similarly to \cite[Lemma 3 and 4]{pavel2019}. 
Then, it follows that $\Phi_p^{-1} \mc A_p$ is $\beta_p \delta_p$-restricted cocoercive with $\delta_p=\frac{1}{|\Phi_p^{-1}|}$ in the $\Phi_p$-induced norm \cite[Lemma 6]{pavel2019}.\\
Concerning $\mc B_p$, it is monotone similarly to \cite[Lemma 4]{pavel2019}.
Consequently, $\Phi_p^{-1} \mc B_p$ is maximally monotone in the $\Phi_p$-induced norm \cite[Lemma 6]{pavel2019}.\\
Since by Assumption \ref{ass_phi_pn}, Assumptions \ref{ass_coco} - \ref{standass_phi} hold, the pFB iterations presented in Algorithm \ref{algo_pe} converge to a v-SGNE of the game in \eqref{eq_game} by Lemma \ref{lemma_FB}.
\qed\end{pf}

\subsection{Convergence of Algorithm \ref{algo_pe}}\label{sec_conv_pe}
Let us now focus on how to obtain Algorithm \ref{algo_pe}. Let us consider the incidence matrix $V$ of the communication graph $\mc G$.
Then, another possibility to force consensus on the dual variables, according to \eqref{eq_null}, is to consider the constraint $\bs V\bs\lambda=0$ (instead of $\bs L\bs \lambda=0$). Exploiting this constraint, we define the two operators 
\begin{equation}\label{eq_op_pe}
\begin{aligned}
\mathcal{C}_p&:\left[\begin{smallmatrix}
\hat{\bs x}\\
\bs v\\
\bs \lambda
\end{smallmatrix}\right]\to\left[\begin{smallmatrix}
\mathcal{R}^{\top} \FF_p(\hat{\bs{x}})+c\bs L_n \hat{\bs{x}}\\
\mathbf{0}_{E m} \\
\bs{b}
\end{smallmatrix}\right]\\
\mc D_p&:\left[\begin{smallmatrix}
\hat{\bs x}\\
\bs v\\
\bs \lambda
\end{smallmatrix}\right]\to\left[\begin{smallmatrix}
\hat G(\hat{\bs{x}}) \\
\mathbf{0}_{\mathrm{Em}} \\
\mathrm{N}_{\mathbb{R}_{\geq 0}^{N m}}(\bs{\lambda})
\end{smallmatrix}\right]+
\left[\begin{smallmatrix}
\mathcal{R}^\top {\bf{A}}^{\top}\lambda\\
-\bs V_m\lambda\\
-{\bf{A}} \mathcal{R}\hat{\bs x}+ \bs V_m \bs v
\end{smallmatrix}\right]
\end{aligned}
\end{equation}
where $\bs{v}=\op{col}((v_{l})_{l \in\{1, \ldots, E\}}) \in \mathbb{R}^{E m}$, $\bs{L}_{n}:=L \otimes I_{n}$ and
$\bs{V}_{m}:=V \otimes \op{Id}_m$. \\
We note that the variable $\bs{v}=\op{col}((v_{l})_{l \in\{1, \ldots, E\}}) \in \mathbb{R}^{E m}$ is used to help reaching consensus on the dual variables. Moreover, it can be interpreted as the network flow. In fact, if we consider $A_i x_i$ as in-flow and $b_i$ as out-flow for each node $i\in\mc I$, then $A\bs x=b$ can be read as a conservative flow balancing constraint. Therefore, $v_l$ can be seen as flow on each edge $l$ to ensure such constraint. In other words, the variable $\bs v$ estimates the contribution of the other agents to the coupling constraints, and ensures that the dual variables reach consensus.\\
Since we consider the incidence matrix $V$ instead of the Laplacian $L$, the preconditioning matrix is given by
\begin{equation}\label{eq_psi_pe}
\Psi_p=\left[\begin{smallmatrix}
\alpha^{-1}  & 0 & -\mc R^\top{\bf{A}}^\top\\
 0 & \nu^{-1} & \bs{V}_m\\
-{\bf {A}}\mc R  & \bs{V}_m^\top & \delta^{-1}
\end{smallmatrix}\right]
\end{equation}
where $\alpha$, $\nu$ and $\delta$ are defined analogously to \eqref{eq_phi_pn}.
Then, by replacing the operator $\mc C_p$ with $\mc C_p^\textup{SA}$, i.e. the operator approximated via the SA scheme according to \eqref{eq_F_SA}, the pFB iteration reads as
\begin{equation}\label{eq_FB_pe}
\bs\omega^{k+1}=(\mathrm{Id}+\Psi^{-1}_p \mc D_p)^{-1} \circ(\mathrm{Id}-\Psi^{-1}_p \mc C_p^\textup{SA})(\bs\omega^k).
\end{equation}
where $\bs\omega^k=\op{col}(\hat{\bs x}^k, \bs v^k,\bs\lambda^k)$.
Then, by expanding \eqref{eq_FB_pe}, 
\begin{equation}\label{eq_FB_pe_1}
\begin{aligned}
\hat{\bs x}^{k+1}&\ni\hat{\bs x}^k -\alpha(\mc R^\top F^\textup{SA}_p(\hat{\bs x})+c\bs L_n\hat{\bs x}+\mc R^\top{\bf{A}}^\top\bs\lambda)+\hat G(\hat{\bs x})\\
v^{k+1}&= v^{k}-\nu\bs V_m \bs{\lambda}^{k} \\
\lambda^{k+1}&= \op{proj}_{\RR_{\geq 0}^{m}}[\lambda^{k}+\delta(2 {\bf{A}}\bs x^{k+1}-{\bf{A}} \bs x^{k}-\bs b\\
&+\bs V^\top_m(2 \bs{v}^{k+1}-\bs{v}^{k}))].
\end{aligned}
\end{equation}
First, we note that also in this case we separate the update of the local decision variables $x_i^k$ and of the estimates $\hat{\bs x}^k_{i,-i}$.
Moreover, in \eqref{eq_FB_pe_1}, two communication rounds are required: one at the beginning of each iteration $k$ to update $\bs x^{k+1}$ and $\bs v^{k+1}$ and one before updating $\bs \lambda^{k+1}$. To avoid this second round, let us introduce, with a little abuse of notation, the variable $\bs z=\op{col}((z_{i})_{i\in\mc I})$ such that $\bs{z}^{k}=\bs V^\top_m \bs{v}^{k},$ for all $k \geq 0 .$ Given an appropriate initialization, e.g., $\bs{v}^{0}=\mathbf{0},$ the following equivalences hold: $\bs{z}^{0}=\mathbf{0}$, $\bs{z}^{k}=\bs V^\top_m \bs{v}^{k}$ and $\bs{z}^{k+1}=\bs V^\top_m \bs{v}^{k+1}=\bs V^\top_m \bs{v}^{k}+\bs V^\top_m \bs V_m \bs{\lambda}^{k}=\bs{z}^{k}+\bs L_m\bs{\lambda}^{k}$. After this change of variables, the iterations in \eqref{eq_FB_pe_1} can be rewritten as in Algorithm \ref{algo_pe}.
Moreover, since the iterations in \eqref{eq_FB_pe} are equivalent to Algorithm \ref{algo_pe}, for the analysis we use the operator in \eqref{eq_op_pe}. Let us also note that the operators are similar to \cite{yi2020} but, due to the change of variable and the approximation scheme, Algorithm \ref{algo_pe} is different from the asynchronous one proposed in \cite{yi2020} and requires less communications. \\
We now proceed in showing that Algorithm \ref{algo_pe} converges to an equilibrium. First, to obtain a positive definite preconditioning matrix (using Gershgorin Theorem), let us bound the step sizes.
\begin{assum}\label{ass_psi_pe}
Given $\tau>0$, the step sizes sequence is such that 
$\bar\alpha$ is as in Assumption \ref{ass_phi_pn} and $\bar \nu$ and $\bar\delta$ are such that, for all $i\in\mc I$
$$0<\nu_i\leq\bar\nu\leq(\tau+\textstyle\sum_{j\in\mc I}\sqrt{w_{ij}})^{-1}$$
$$0<\delta_{i} \leq\bar\delta\leq(\tau+\textstyle\sum_{j\in\mc I}\sqrt{w_{ij}}+\max _{j\in\{1, \ldots m\}}\textstyle\sum_{k=1}^{n_{i}}|\left[A_{i}\right]_{j k}|)^{-1}$$
where $[A_i^\top]_{jk}$ indicates the entry $(j,k)$ of the matrix $A_i^\top$ and such that $\Psi_p$ satisfies Assumption \ref{standass_phi}.
\end{assum}

Finally, we prove Theorem \ref{theo_pe}.\vspace{-.5cm}
\begin{pf}[Proof of Theorem \ref{theo_pe}]
First, we relate the unique v-SGNE of the game in \eqref{eq_game} to the zeros of $\mathcal{C}_p+\mc D_p$ in \eqref{eq_op_pe}.
Namely, given any $\hat{\bs \omega}^{*}:=\op{col}(\hat{\bs x}^{*}, \bs{v}^{*}, \bs \lambda^{*}) \in \op{zer}(\mathcal{C}_p+\mc D_p),$ it holds that $\hat{\bs x}^{*}=\mathbf{1}_{N} \otimes \bs x^{*}$ and $\bs \lambda^{*}=\mathbf{1}_{N} \otimes \lambda^{*}$, where the pair $(\bs x^{*}, \lambda^{*})$ satisfies the KKT conditions \eqref{eq_T}, i.e., $\bs x^{*}$ is a v-SGNE of the game in \eqref{eq_game}.
This follows analogously to Theorem \ref{theo_pn} noting that, expanding $\hat{\bs\omega}\in\op{zer}(\mc C_p+\mc D_p)$, from the second line we have $\bs \lambda\in \op{null}(\bs V)$ and we can premultiply by $(\bs 1_E^\top\otimes\op{Id}_m)$ the third line to obtain the KKT conditions in \eqref{eq_T}. 
Similarly, $\op{zer}(\mathcal{C}_p+\mc D_p) \neq \varnothing$.\\
Then, we show the monotonicity properties of the operators.
Note that $\mc C_p$ is the same as $\mc A_p$ therefore it is $\beta_p$-restricted cocoercive with respect to $Z$, where $\beta_p\in(0,\frac{\mu_p}{\theta_p^2}]$, $\mu=\uplambda_{\min}(\Upsilon_p)$ and $\theta_p=\ell_{p}+2cd_{\text{max}}$ by Theorem \ref{theo_pn}. Therefore, $\Psi_p^{-1} \mc C_p$ is $\beta_p \delta_p$-restricted cocoercive with respect to $Z$, with $\delta_p>\frac{1}{|\Psi_p^{-1}|}$ in the $\Psi_p$-induced norm.\\
The operator $\mc D_p$ is maximally monotone analogously to the proof of Theorem \ref{theo_pn}. It follows that $\Psi_p^{-1} \mc D_p$ is maximally monotone in the $\Psi_p$-induced norm.\\
Convergence follows by Lemma \ref{lemma_FB}.
\qed\end{pf}

\section{Convergence analysis for aggregative games}\label{sec_conv_a}
Analogously to Section \ref{sec_conv_p}, we show in this section that the algorithms proposed for the aggregative case in Section \ref{sec_agg} converge a.s. to a v-SGNE of the game in \eqref{eq_game}.\\
We recall that we keep track of the aggregate value through the variable $\bs s=1_N\otimes \op{avg}(\bs x)-\bs x$ and that the approximated average is given by $\bs u=\bs x+\bs s$ as in \eqref{eq_u}.

\subsection{Convergence of Algorithm \ref{algo_an}}\label{sec_conv_an}
Let us start by defining the operators that leads to the iterations in Algorithm \ref{algo_an}. Specifically, the forward and backward operators should be defined according to
\begin{equation}\label{eq_op_an}
\begin{aligned}
\mc A_a&:\left[\begin{smallmatrix}
\bs{x} \\
\bs s\\
\bs z\\
\bs \lambda
\end{smallmatrix}\right] \mapsto\left[\begin{smallmatrix}
\FF_a(\bs x, \bs x+\bs s)+c\bs{L}_{\bar{n}} (\bs x+\bs s) \\
\bs{L}_{\bar{n}} (\bs x+\bs s)\\
\mathbf{0}_{N m}\\
\bs{b}
\end{smallmatrix}\right]\\
\mc B_a&:\left[\begin{smallmatrix}
\bs{x} \\
\bs s\\
\bs z\\
\bs \lambda
\end{smallmatrix}\right] \mapsto\left[\begin{smallmatrix}
G(\bs x) \\
\mathbf{0}_{n} \\
\mathbf{0}_{N m }\\
\mathrm{N}_{\mathbb{R}_{\geq 0}^{Nm}}(\bs{\lambda})
\end{smallmatrix}\right]+
\left[\begin{smallmatrix}
{\bf{A}}^{\top}\lambda\\
\bs 0\\
-\bs L_m\lambda\\
-{\bf{A}} \hat{\bs x}+ \bs L_m \bs z
\end{smallmatrix}\right]
\end{aligned}
\end{equation}
where $\bs L_{\bar n}=L\otimes \op{Id}_{\bar n}$, $\bs L_m=L\otimes \op{Id}_m$.
We note that, compared to the operators $\mc A_p$ and $\mc B_p$ in \eqref{eq_op_pn}, in \eqref{eq_op_an} instead of the estimates we can take the true decision variables $\bs x$ as a state variable because we track the average value through the variable $\bs s$.
The preconditioning matrix reads similarly to $\Phi_p$ in \eqref{eq_phi_pn} with the addition of a line (corresponding to the variable $\bs s$):
\begin{equation}\label{eq_phi_an}
\Phi_a=\left[\begin{smallmatrix}
\alpha^{-1} & 0 & 0 & -{\bf{A}}^\top\\
0 & \gamma^{-1} & 0 & 0\\
0 & 0 & \nu^{-1} & \bs{L}_m\\
-{\bf{A}} & 0 & \bs{L}_m & \delta^{-1}
\end{smallmatrix}\right],
\end{equation} 
where $\alpha^{-1}=\op{diag}\{\alpha_1^{-1}\op{Id}_{\bar n},\dots,\alpha_N^{-1}\op{Id}_{\bar n}\}\in\RR^{n\times n}$ and similarly $\gamma^{-1}$, $\nu^{-1}$ and $\delta^{-1}$ are block diagonal matrices of suitable dimensions.
Given the operators in \eqref{eq_op_an} and the preconditioning matrix in \eqref{eq_phi_an}, Algorithm \ref{algo_an} in compact form reads as the pFB iteration
\begin{equation*}\label{eq_fb_an}
\bs\omega^{k+1}=(\mathrm{Id}+\Phi_a^{-1} \mc B_a)^{-1} \circ(\mathrm{Id}-\Phi_a^{-1} \mc A^{\textup{SA}}_a)(\bs\omega^k,\xi^k),
\end{equation*}
where $\bs\omega^k=\op{col}(\bs x^k,\bs s^k,\bs z^k,\bs\lambda^k)$, $(\mathrm{Id}+\Phi^{-1} \mc B_a)^{-1}$ represent the backward step and $(\mathrm{Id}-\Phi^{-1}\mc A^{\textup{SA}}_a)$ is the forward step where the pseudogradient mapping is approximated according to \eqref{eq_F_SA}.
\begin{remark}
In line with Remark \ref{remark_gadjov}, let us note that in \cite{gadjov2020} a different splitting and preconditioning are used. Since the paper tracks the aggregate value $\bs u$ (instead of $\bs s$), to prove convergence, it uses an auxiliary iterative scheme on the orthogonal complement of the consensus subspace \cite[Lemma 2]{gadjov2020}. Hence, the extended operators and preconditioning matrix depend on the projection on the consensus subspace ($P_\parallel$) and on its orthogonal complement ($P_\perp$), while ours are a generalization of the operators in \eqref{eq_op_pn} and \eqref{eq_phi_pn} to the aggregative case. We can avoid the auxiliary iteration because we track the aggregative value with $\bs s$ instead of measuring the disagreement with $\bs u^\perp$. We also note that, in general, $\bs s\neq\bs u^\perp$.
\end{remark}
Now we prove that using the operators $\mc A_a$ and $\mc B_a$ in \eqref{eq_op_an} we can reach a v-SGNE of the game in \eqref{eq_game}. To this aim, we restrict our analysis to the invariant subspace
\begin{equation}\label{eq_inv_an}
\Sigma:=\{(\bs x, \bs s, \bs{z}, \bs{\lambda}) \in \mathbb{R}^{2n+2Nm} | \op{avg}(\bs s)=\mathbf{0}_{\bar{n}}\},
\end{equation}
i.e., the space where the agents are able to compute the exact aggregative value.
Moreover, to ensure that $\mc A^a$ and $\mc B^a$ have the monotonicity properties of Lemma \ref{lemma_FB}, the preconditioning matrix should be positive definite. For this reason, we take some bounds on the step sizes.
\begin{assum}\label{ass_phi_an}
The step sizes $\bar\alpha$, $\bar\nu$ and $\bar\delta$ satisfy Assumption \ref{ass_phi_pn} and $\gamma>0$. Moreover, they are such that $\Phi_a$ satisfies Assumption \ref{standass_phi}.
\end{assum}

We are ready to prove convergence. \vspace{-.5cm}
\begin{pf}[Proof of Theorem \ref{theo_an}]
First, we ensure that the zeros of $(\mc A_a+\mc B_a)\cap\Sigma$ are v-SGNEs. 
Let us consider $\bs\omega^*=\op{col}(\bs x^*,\bs s^*,\bs z^*,\bs\lambda^*)\in\op{zer}(\mc A_a+\mc B_a)\cap\Sigma$ and for brevity, $\bs u^*=\bs x^*+\bs s^*$, then
\begin{equation}\label{eq_zeros_an}
\begin{aligned}
0 & \in \FF_a(\bs x^{*}, \bs u^{*})+c\bs L_{\bar{n}} \bs u^{*}+G(\bs x)+{\bf{A}}^{\top} \boldsymbol{\lambda}^{*} \\
0 &=\bs L_{\bar{n}} \bs u^{*} \\
0&=\boldsymbol{L}_{m} \boldsymbol{\lambda}^{*} \\
0 & \in \boldsymbol{b}+\mathrm{N}_{\mathbb{R}_{\geq 0}^{N m}}(\boldsymbol{\lambda}^{*})-{\bf{A}} x^{*}+\boldsymbol{L}_{m}^{\top} \boldsymbol{z}^{*}
\end{aligned}
\end{equation}
Let us recall that ${\bf{A}}^{\top}(\mathbf{1}_{N} \otimes \lambda^{*})=A^{\top} \lambda^{*}$ and $\FF_a(\bs x^{*}, \bs{1}_{N} \otimes \bs x^{*})=\FF(\bs x^{*})$, then, from the first line of \eqref{eq_zeros_an}, we obtain the first line of the KKT conditions in \eqref{eq_T}.
From the third line of \eqref{eq_zeros_an} and from \eqref{eq_null}, we have $\bs \lambda^{*}=1_{N} \otimes \lambda^{*},$ for some $\lambda^{*} \in \mathbb{R}^{m}_{\geq 0}$. 
From the second line and since $\bs\omega^{*} \in \Sigma,$ it holds that $\bs u^{*}=\bs x^{*}+\bs s^{*}=\bs 1_{N} \otimes \op{avg}(\bs x^{*})$. 
Since $(\mathbf{1}_{N}^{\top} \otimes \op{Id}_m) \bs{b}=b$, $(\mathbf{1}_{N}^{\top} \otimes \op{Id}_m) \bs{L}_{m}=0$ (by \eqref{eq_null} and symmetry of L), $(\mathbf{1}_{N} \otimes \op{Id}_m) {\bf{A}}=A$ and $(\mathbf{1}_{N}^{\top} \otimes \op{Id}_m) N_{\mathbb{R}_{>0}^{N m}}(\mathbf{1}_{N} \otimes \lambda^{*})=N \mathrm{N}_{\mathbb{R}_{\geq0}^{m}}(\lambda^{*})=\mathrm{N}_{\mathbb{R}_{\geq0}^{m}}(\lambda^{*})$,
we premultiply the fourth line by $(\mathbf{1}_{N}^{\top} \otimes \op{Id}_m)$ to obtain the second line of \eqref{eq_T}.
Therefore, the pair $(\bs x^{*}, \lambda^{*})$ satisfies the KKT conditions \eqref{eq_T}, i.e., $\bs x^{*}$ is a v-SGNE of the game in \eqref{eq_game}.\\
Moreover, $\op{zer}(\mc A_a+\mc B_a) \cap \Sigma \neq \varnothing$.
From Assumption \ref{ass_strong}, it follows that there is only one v-SGNE, i.e., a pair $(\bs x^*,\bs \lambda^*)$ that satisfy the KKT conditions in \eqref{eq_T}. Now, we show that there exists $\bs{z}^{*} \in \mathbb{R}^{N m}$ such that $\bs{\omega}^{*}=$ $\op{col}(\bs x^{*}, \mathbf{1}_{N} \otimes \op{avg}(\bs x^{*})-\bs x^{*}, \bs{z}^{*}, \mathbf{1}_{N} \otimes \lambda^{*}) \in \op{zer}(\mc A_a+\mc B_a) \cap \Sigma$. 
It holds that $\bs\omega^{*} \in \Sigma$ and that $\bs\omega^{*}$ satisfies the first three lines of \eqref{eq_zeros_an}. Exploiting the KKT conditions in \eqref{eq_T}, there exists $w^{*} \in \mathrm{N}_{\mathbb{R}_{>0}^{m}}(\lambda^{*})$ such that $A x^{*}-b-w^{*}=\mathbf{0}_{n} .$ Moreover, $\mathrm{N}_{\mathbb{R}_{\geq 0}^{N m}}(1_{N} \otimes \lambda^{*})=\Pi_{i \in \mathcal{I}} \mathrm{N}_{\mathbb{R}_{\geq 0}^{m}}(\lambda^{*})$
and it follows by properties of the normal cone that $\op{col}(w_{1}^{*}, \ldots, w_{N}^{*}) \in\mathrm{N}_{\mathbb{R}_{\geq 0}^{N m}}(\mathbf{1}_{N} \otimes \lambda^{*})$, with $w_{1}^{*}=\cdots=w_{N}^{*}=\frac{1}{N} w^{*}$. Therefore $(\mathbf{1}_{N}^{\top} \otimes \op{Id}_m)(-{\bf{A}} x^{*}+\bs{b}+\op{col}(w_{1}^{*}, \ldots, w_{N}^{*}))=b-A x^{*}+w^{*}=$
$\mathbf{0}_{m},$ or $-{\bf{A}} x^{*}+\bs{b}+\op{col}(w_{1}^{*}, \ldots, w_{N}^{*}) \in \op{null}(\mathbf{1}_{N}^{\top} \otimes \op{Id}_m) \subseteq\op{range}(L_m)$. Since $ \op{range}(\bs{L}_{m})=\op{null}(\mathbf{1}_{N}^{\top} \otimes I_{m})=C_{m}^{\perp}$, there exists $\bs{z}^{*}$ such that also the last line of \eqref{eq_zeros_an} is satisfied, i.e., $\bs\omega^{*} \in \op{zer}(\mc A_a+\mc B_a)\cap\Sigma$.
We now prove the monotonicity properties of the operators $\mc A_a$ and $\mc B_a$.\\
Similarly to the proof of Theorem \ref{theo_pn} and \cite[Lemma 4]{gadjov2020}, let us introduce
$$\Upsilon_a=\left[\begin{smallmatrix}
\eta & -\frac{\ell^u_a}{2}\\
-\frac{\ell^u_a}{2} & \uplambda_2(L)
\end{smallmatrix}\right],$$
where $\eta$ is the strong monotonicity constant as in Assumption \ref{ass_strong} and $\ell^u_a$ the Lipschitz constant as in Remark \ref{remark_lip_a}.
Then, the two operators $\mc A_a$ and $\mc B_a$ in \eqref{eq_op_an} have the properties of Assumption \ref{ass_coco}. To prove this
let us note that each vector can be decomposed as $\bs u=\bs u^\parallel+\bs u^\perp$ where $\bs u^\parallel=P_\parallel\bs u\in C_{\bar n}$ and $\bs u^\perp=P_\perp\bs u\in C_{\bar n}^\perp$ and $P_\parallel$ and $P_\perp$ are the projection operators defined as $P_\parallel=\frac{1}{N}\bs 1_N\otimes\bs 1_N^\top\otimes \op{Id}_{\bar n}$ and $P_\perp=\op{Id}_{Nn}-P_\parallel$ \cite{gadjov2020}. We note that $P_\parallel \bs u^k=\bs 1_N\otimes \op{avg}(\bs u^k)=P_\parallel\bs x^k$ by the property in \eqref{eq_invariance}.\\
First, we prove that $\mc A_a$ is $\beta_a$-restricted cocoercive with constant $0<\beta_a \leq\frac{\mu_a}{\theta_a}$, where $\mu_a=\uplambda_{\min}(\Upsilon_a)$ and $\theta_a=\max\{(2\ell_a^x)^2,(2\ell_a^u)^2+3\uplambda_2(L)\}$. To this aim, we first prove that the operator
\begin{equation*}
\bar{\mc A}_a:\left[\begin{smallmatrix}
\bs x\\
\bs s
\end{smallmatrix}\right]\mapsto
\left[\begin{smallmatrix}
\FF(\bs x, \bs u)+ \bs{L}_{\bar{n}} \bs u\\
\bs{L}_{\bar{n}} \bs u \\
\end{smallmatrix}\right]
\end{equation*}
is restricted strongly monotone, then restricted cocoercivity of $\mc A_a$ follows.
Let $\bs w=\op{col}(\bs x, \bs s)$, then
\begin{equation*}
\begin{aligned}
&\langle \bar{\mc A}_a\bs w-\bar{\mc A}_a\bs w^*,\bs w-\bs w^*\rangle=\\
&=\langle \FF_a(\bs x,\bs u)-\FF_a(\bs x^*,\bs u^*),\bs x-\bs x^*\rangle+\langle L\bs u-L\bs u^*,\bs u-\bs u^*\rangle\\
&+\langle \FF_a(\bs x,\op{avg}(\bs x))-\FF_a(\bs x,\op{avg}(\bs x)),\bs x-\bs x^*\rangle
\end{aligned}
\end{equation*}
using the fact that $\bs u=\bs x+\bs s$.
Moreover, $\langle \bs L_{\bar n}\bs u-\bs L_{\bar n}\bs u^*,\bs u-\bs u^*\rangle\geq\uplambda_2(L)\normsq{\bs u^\perp}$ as in \cite[Lemma 4]{gadjov2020}.
Notice that $\bs u-\op{avg}(\bs x)=\bs u^\parallel+\bs u^\perp-\op{avg}(\bs x)=\bs u^\perp$ by the invariance property \eqref{eq_invariance}.
Therefore, similarly to \cite[Lemma 4]{gadjov2020},
\begin{equation}\label{step1}
\begin{aligned}
\langle \bar{\mc A}_a\bs w&-\bar{\mc A}_a\bs w^*,\bs w-\bs w^*\rangle\geq \\
&\geq\left[\begin{smallmatrix}
\norm{\bs x-\bs x^*}\\
\norm{\bs u^\perp}
\end{smallmatrix}\right]^\top\Upsilon_a\left[\begin{smallmatrix}
\norm{\bs x-\bs x^*}\\
\norm{\bs u^\perp}
\end{smallmatrix}\right].\\
\end{aligned}
\end{equation}
To prove cocoercivity, we use an argument similar to \cite[Lemma 4]{gadjov2020} to obtain
\begin{equation}\label{step2}
\normsq{\bar{\mc A}_a\bs w-\bar{\mc A}_a\bs w^*}\leq\theta_a\left\|\left[\begin{smallmatrix}
\bs x-\bs x^*\\
\bs u^\perp
\end{smallmatrix}\right]\right\|^2
\end{equation}
The constant $\theta_a$ is defined as $\theta_a=\max\{(2\ell_a^x)^2,(2\ell_a^u)^2+3\uplambda_2(L)\}$. Pairing \eqref{step1} and \eqref{step2}, we have that $\bar{\mc A}_a$ is $\frac{\mu_a}{\theta_a}$-cocoercive with $\mu_a=\uplambda_{\min}(\Upsilon_a)$.
The fact that $\Phi^{-1}\mc A_a$ is restricted cocoercive follows from $\mc A_a$ being restricted cocoercive and the maximal monotonicity of $\mc B_a$ and $\Phi^{-1}\mc B_a$ follows analogously to Theorem \ref{theo_pn}. Since Assumptions \ref{ass_coco}--\ref{standass_phi} are verified, convergence follows from Lemma \ref{lemma_FB}.
\qed\end{pf}

\subsection{Convergence of Algorithm \ref{algo_ae}}\label{sec_conv_ae}
We now consider Algorithm \ref{algo_ae} and show how to obtain its iterates. Later we prove Theorem \ref{theo_ae}.\\
Including the variable $\bs s$ to track the aggregative value and considering the consensus constraint $\bs V\bs\lambda=0$, the operators of the edge-based algorithm read as
\begin{equation}\label{eq_op_ae}
\begin{aligned}
\mc C_a&:\left[\begin{smallmatrix}
\bs{x} \\
\bs s\\
\bs v\\
\bs \lambda
\end{smallmatrix}\right]\to\left[\begin{smallmatrix}
\FF_a(\bs x, \bs x+\bs s)+c\bs{L}_{\bar{n}} (\bs x+\bs s)\\
\bs{L}_{\bar{n}} (\bs x+\bs s) \\
\mathbf{0}_{E m} \\
\bs{b}
\end{smallmatrix}\right]\\
\mc D_a&:\left[\begin{smallmatrix}
\bs{x} \\
\bs s\\
\bs v\\
\bs \lambda
\end{smallmatrix}\right]\to\left[\begin{smallmatrix}
G(\bs x) \\
\mathbf{0}_{n} \\
\mathbf{0}_{E m} \\
\mathrm{N}_{\mathbb{R}_{>0}^{Nm}(\bs{\lambda})}
\end{smallmatrix}\right]+\left[\begin{smallmatrix}
{\bf{A}}^{\top} \bs{\lambda} \\
\mathbf{0}_{\bs{V}} \\
-\bs{V}_{m} \bs{\lambda} \\
\bs V^\top_m\bs{v}-{\bf{A}} \bs x
\end{smallmatrix}\right],
\end{aligned}
\end{equation}
while the preconditioning matrix is given by
\begin{equation}\label{eq_psi_ae}
\Psi_a=\left[\begin{smallmatrix}
\alpha^{-1} & 0 & 0 & -{\bf{A}}^\top\\
0 & \gamma^{-1} & 0 & 0\\
0 & 0 & \nu^{-1} & -\bs{V}_m\\
-{\bf{A}} & 0 & -\bs{V}_m & \delta^{-1}
\end{smallmatrix}\right],
\end{equation} 
where $\bs V_m=V\otimes\op{Id}_m$ and $\alpha$, $\gamma$, $\nu$ and $\delta$ are the block diagonal step sizes. 
Similarly to the previous sections, the pFB iteration is given by
\begin{equation}\label{eq_FB_ae}
\bs\omega^{k+1}=(\mathrm{Id}+\Psi^{-1}_a \mc D_a)^{-1} \circ(\mathrm{Id}-\Psi^{-1}_a \mc C_a^\textup{SA})(\bs\omega^k),
\end{equation}
where $\mc C_a^\textup{SA}$ is approximated according to \eqref{eq_F_SA}. The iterations in Algorithm \ref{algo_ae} can be obtained expanding \eqref{eq_FB_ae} and using a change of variables $\bs{z}^{k}=\bs V^\top_m \bs{v}^{k},$ as in Section \ref{sec_pe}.\\
To prove convergence of Algorithm \ref{algo_ae}, we consider the invariant subspace
$$
\Sigma:=\{(\bs x, \bs s, \bs{v}, \bs{\lambda}) \in \mathbb{R}^{2n+Em+Nm}| \op{avg}(\bs s)=\mathbf{0}_{\bar{n}}\},
$$
similarly to Section \ref{sec_an}. We note that the dimension of this set is different from \eqref{eq_inv_an}.
Before stating the result, analogously to the previous sections, we must guarantee that the preconditioning matrix is positive definite.
\begin{assum}\label{ass_psi_ae}
Given $\tau>0$, the step sizes are such that $\gamma>0$,
$\bar\alpha$ is as in Assumption \ref{ass_phi_pn}, $\bar\nu$ and $\bar\delta$ are such that, for all $i\in\mc I$,
$$0<\nu_i\leq\bar\nu\leq(\tau+\textstyle\sum_{j\in\mc I}\sqrt{w_{ij}})^{-1},$$
$$0<\delta_{i} \leq\bar\delta\leq(\tau+\textstyle\sum_{j\in\mc I}\sqrt{w_{ij}}+\max _{j\in\{1, \ldots m\}}\textstyle\sum_{k=1}^{n_{i}}|\left[A_{i}\right]_{j k}|)^{-1},$$
where $[A_i^\top]_{jk}$ indicates the entry $(j,k)$ of the matrix $A_i^\top$ and such that $\Psi_a$ satisfies Assumption \ref{standass_phi}.
\end{assum}

Then, the convergence result holds.

\begin{pf}[Proof of Theorem \ref{theo_ae}]
We start by showing that, given any $\bs\omega=\op{col}(\bs x^{*}, \bs s^{*}, \bs{v}^{*}, \bs{\lambda}^{*}) \in \op{zer}(\mc C_a+\mc D_a) \cap \Sigma,$ it holds that $\bs s^{*}=\mathbf{1}_{N} \otimes\op{avg}(\bs x^{*})-\bs x^{*}$, $\bs \lambda^{*}=\mathbf{1}_{N} \otimes \lambda^{*}$ and the pair $(\bs x^{*}, \lambda^{*})$ satisfies the KKT conditions \eqref{eq_T}, i.e., $\bs x^{*}$ is a v-SGNE of the game in \eqref{eq_game}.\\
To this aim, let us consider any $\bs{\omega}^{*} \in \op{zer}(\mc C_a+\mc D_a) \cap \Sigma$ and let $\bs u^{*}=\bs x^{*}+\bs s^{*} ;$ then we have
\begin{equation*}\label{eq_zeros}
\begin{aligned}
&\mathbf{0}_{\bar{n}} \in  \FF_a(\bs x^{*}, \bs u^{*})+\bs{L}_{\bar{n}} \bs u^{*}+G(\bs x^{*})+{\bf{A}}^{\top} \bs{\lambda}^{*} \\
&\mathbf{0}_{\bar{n}} =\bs{L}_{\bar{n}} \bs u^{*} \\
&\mathbf{0}_{E m} =\bs{V}_{m} \bs{\lambda}^{*} \\
&\bs{0}_{N m} \in \bs{b}+\mathrm{N}_{\mathbb{R}_{\geq0}^{N m}}(\bs{\lambda}^{*})-{\bf{A}} \bs x^{*}-\bs{V}_{m}^{\top} \bs{v}^{*}
\end{aligned}
\end{equation*}
The fact that the KKT conditions in \eqref{eq_T} are satisfied, i.e., that $\bs x^*$ is a v-SGNE of the game in \eqref{eq_game}, follows analogously to Theorem \ref{theo_an}.\\
Moreover, it holds that $\op{zer}(\mc C_a+\mc D_a) \cap \Sigma \neq \varnothing$, similarly to Theorem \ref{theo_an} and using the fact that $\op{range}(\bs{V}_{q}^{\top}) \supseteq \op{range}(\bs{L}_{q})=\op{null}(\mathbf{1}_{N}^{\top} \otimes I_{q})=\bs{E}_{q}^{\perp}$.\\
The fact that
$\mc C_a$ is $\beta_a$-restricted cocoercive where $0<\beta_a \leq\frac{\mu_a}{\theta_a}$, $\mu_a=\uplambda_{\min}(\Upsilon_a)$ and $\theta_a=\max\{(2\ell_a^x)^2,(2\ell_a^u)^2+3\uplambda_2(L)\}$ and that the operator $\mc D_a$ is maximally monotone and that $\Psi_a^{-1}\mc C_a$ is $\beta_a\delta_a$-cocoercive with $\delta_a=\frac{1}{|\Phi_a^{-1}|}$ and that
$\Psi_a^{-1}\mc D_a$ is maximally monotone in the $\Phi_a$-induced norm follows analogously to Theorem \ref{theo_an}. Convergence follows by Lemma \ref{lemma_FB}.
\qed\end{pf}

\section{Numerical simulations}\label{sec_sim}
Let us present some numerical results to validate the convergence analysis of our proposed algorithms. For network games, we consider a Nash-Cournot game \cite{koshal2016,pavel2019,yu2017} while for the aggregative case we use a charging scheduling problem \cite{grammatico2017}. In both cases we consider two different topologies for the communication graph, i.e., a complete graph and a cycle graph, to show how connectivity affects the results.
For the tests, we take the step sizes to be half of those that generate instability.

\subsection{Nash-Cournot games}
\begin{figure}[t]
\centering
\includegraphics[width=\columnwidth]{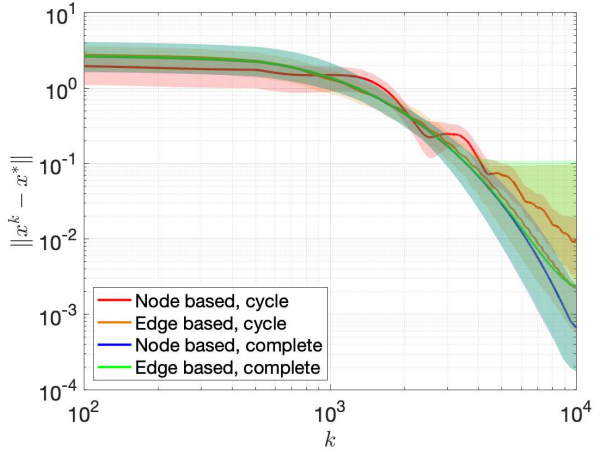}
\caption{Distance from a solution (network game).}\label{fig_net_iter}
\end{figure}
\begin{figure}[t]
\centering
\includegraphics[width=\columnwidth]{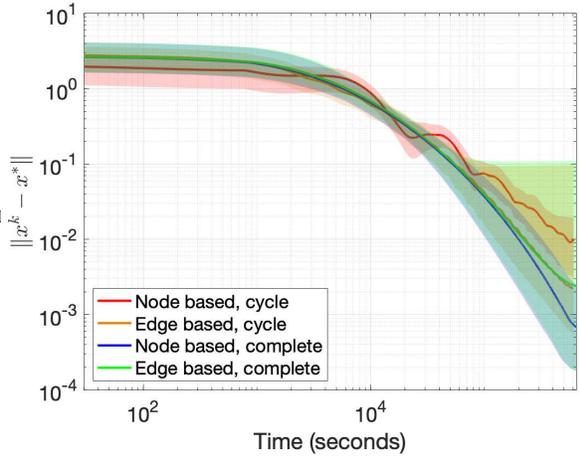}
\caption{Computational time (network game).}\label{fig_net_tempo}
\end{figure}
In a Nash-Cournot game \cite{koshal2016,pavel2019,yu2017} a set $\mc I=\{1,\dots,N\}$ of $N=20$ companies (agents) produce a commodity to be sold over $m=7$ markets as in \cite[Figure 1]{pavel2019}. The markets have a bounded capacity therefore the companies face some coupling constraints.
We suppose each company $i$ to have a strongly convex, quadratic cost of production $c_{i}(x_{i},\xi_i)=x_{i}^{T} Q_{i}(\xi_i) x_{i}+q_{i}^{T} x_{i}$ where $Q_{i}(\xi_i)$ is a random diagonal matrix with the entries drawn from a normal distribution with mean $4.5$ and bounded variance, while each component of $q_{i}\in\RR^{n_i}$ is taken from $[1,2]$.
Each market $j$ has a linear price depending on the total amount of commodities sold to it: $P_j(x,\zeta)=p_{j}(\zeta)-\chi_{j}[A x]_{j}$ with $p_j(\zeta)$ a random vector drawn from a normal distribution with mean 15 and bounded variance and $\chi_j\in[1,3]$. The cost function of each agent reads as
$\JJ_i(x_i,\bs x_{-i})=\EE[c_i(x_i,\xi_i)-P(\bs x,\xi)^\top A_ix_i],$
and it satisfies the assumptions of our problem since it is strongly convex \cite[Section V-A]{pavel2019}. 
We suppose that the local constraint of the companies are given by some bounds on the production, i.e., $\Omega_i=\{x_i\in\RR^{n_i}:0\leq x_i\leq X_i\}$, $i\in\mc I$, where each $X_i$ is randomly drawn in $[5,10]$. Each market $j$ has a maximal capacity of $b_j,$ randomly drawn from $[1,2]$. 
In Figure \ref{fig_net_iter} and \ref{fig_net_tempo} we show the results for network games. Specifically, Figure \ref{fig_net_iter} we show the distance from the solution versus the number of iterations while Figure \ref{fig_net_tempo} shows the computational cost; the transparent areas show the variance over 100 runs of the algorithm. We discard the first 100 iterations to better visualize the asymptotic behavior. As one can see, there is not a significant difference between the two algorithms for the complete case while the node-based algorithm is slower for the cycle graph.

\subsection{Charging scheduling problem}

\begin{figure}[t]
\centering
\includegraphics[width=\columnwidth]{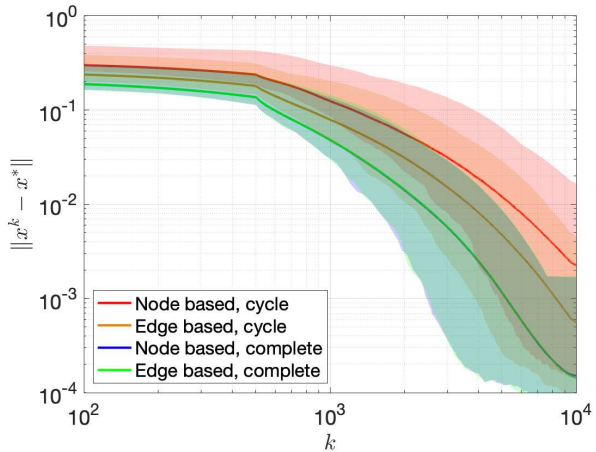}
\caption{Distance from a solution (aggregative game).}\label{fig_agg_iter}
\end{figure}

\begin{figure}[t]
\centering
\includegraphics[width=\columnwidth]{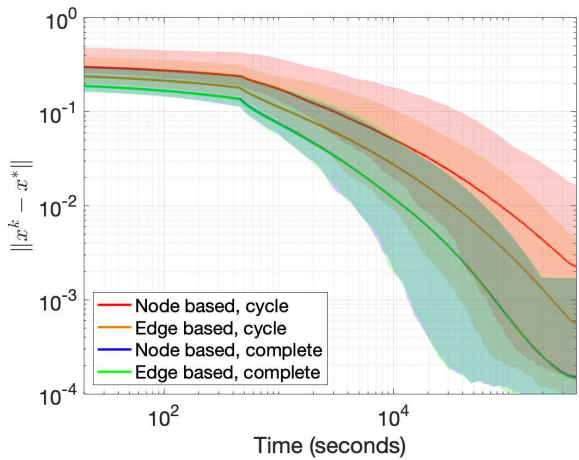}
\caption{Computational time (aggregative game).}\label{fig_agg_tempo}
\end{figure}

For aggregative games, we consider the charging scheduling problem, where the agents are plug-in vehicles inspired by \cite{grammatico2017,lei2020}. We suppose to have $N=10$ users, planning the charging profile over the next 24 hours, divided into $\bar n=12$ time slots. Each user has a random linear battery degradation cost $c_{i}(x_{i}, \xi_{i})=c_i(\xi_i)^\top x_i$ for some random vector $c_i(\xi_i)$ drawn from a normal distribution with mean 4 and bounded variance.
The cost of energy for each time slot depends on the aggregate value, i.e., $P_j(\op{avg}(\bs x),\zeta)=p_j(\zeta)-\chi_j[\op{avg}(\bs x)]_j$. The random variables $p_j(\zeta)$ are drawn from a normal distribution with mean 4.5 and bounded variance while $\chi_j\in[1,2]$.
Therefore each agent has a cost function of the form
$\JJ_i(x_i,\bs x_{-i})=\EE[c_i(x_i,\xi_i)-P(\op{avg}(\bs x),\zeta)^\top x_i].$
The local constraints are given by $\Omega_i=\{x_i\in\RR^{n_i}:0\leq x_i\leq X_i\}$ where $X_i$ is taken according to the following rule: $[X_i]_j=0.25$ with probability $1/2$ and $[X_i]_j=0$ otherwise. Moreover, the users are subject to the transition line constraint $0\leq\sum_{i\in\mc I}x_i\leq\theta$ where $\theta_j=0.4$ if $j\in\{1,2,3,11,12\}$, i.e., it is less restrictive at night, and $\theta_j=1$ during the day.
In Figures \ref{fig_agg_iter} and \ref{fig_agg_tempo}, we show the distance from a solution; the transparent areas show the variance over 100 runs of the algorithm. We discard the first 100 iterations to better visualize the asymptotic behavior. From these figures we see that also in this case the node-base algorithm for the cycle graph is the slowest in terms of both number of iterations (Figure \ref{fig_agg_iter}) and computational time (Figure \ref{fig_agg_tempo}).



\section{Conclusion}\label{sec_conclu}
The preconditioned forward-backward (pFB) algorithm can be used to find stochastic generalized Nash equilibria in a partial-decision information setup. Leveraging on the estimation of the unknown variables, the pFB algorithm can be tailored for network games and for aggregative games. Thanks to the preconditioning almost sure convergence holds under restricted cocoercivity of the forward operator with respect to the solution set.

\balance

\bibliographystyle{abbrv}        
\bibliography{biblio}           

\end{document}